\documentclass[a4paper,reqno,11pt]{amsart}

\usepackage{fullpage}

\usepackage{amssymb}
\usepackage{dsfont}
\usepackage{mathrsfs}

\usepackage[utf8]{inputenc} 
\usepackage[T1]{fontenc} 

\usepackage{hyperref}
\usepackage[nobysame,alphabetic,initials,msc-links]{amsrefs}

\DefineSimpleKey{bib}{how}

\renewcommand{\eprint}[1]{#1}
\BibSpec{misc}{%
  +{}{\PrintAuthors}  {author}
  +{,}{ \textit}      {title}
  +{,}{ }             {how}
  +{}{ \parenthesize} {date}
  +{,} { available at \eprint}        {eprint}
  +{,}{ available at \url}{url}
  +{,}{ }             {note}
  +{.}{}              {transition}
}

\numberwithin{equation}{section}

\theoremstyle{plain}
\newtheorem{thm}{Theorem}[section]
\newtheorem{prop}[thm]{Proposition}
\newtheorem{lemma}[thm]{Lemma}
\newtheorem{cor}[thm]{Corollary}
\theoremstyle{definition}

\newtheorem{defn}[thm]{Definition}
\theoremstyle{remark}

\newtheorem{remark}[thm]{Remark}
\newtheorem{example}[thm]{Example}

\theoremstyle{plain}

\newcommand\bp{\begin{proof}}
\newcommand\ep{\end{proof}}

\newcommand{\un}{\mathds{1}}

\newcommand\C{\mathbb{C}}
\newcommand\N{\mathbb{N}}

\newcommand\R{\mathbb{R}}
\newcommand\T{\mathbb{T}}
\newcommand\Z{\mathbb{Z}}

\newcommand{\A}{\mathcal{A}}
\newcommand{\BB}{\mathcal{B}}

\newcommand{\G}{\mathcal{G}}
\newcommand{\HH}{\mathcal{H}}
\newcommand{\LL}{\mathcal{L}}
\newcommand{\M}{\mathcal{M}}
\newcommand{\NN}{\mathcal{N}}

\newcommand{\UU}{\mathcal{U}}

\newcommand\Ad{\operatorname{Ad}}
\newcommand\Bis{\operatorname{Bis}}
\newcommand\BisG{\operatorname{Bis}(\mathcal G)}
\newcommand\BiscG{\operatorname{Bis}_c(\mathcal G)}
\newcommand\BisGL{\operatorname{Bis}(\mathcal G,\mathcal L)}

\newcommand{\Gu}{{\mathcal{G}^{(0)}}}
\newcommand{\Gxx}{\mathcal{G}^x_x}

\newcommand\supp{\operatorname{supp}}
\newcommand\Ped{\operatorname{Ped}}

\newcommand\eps{\varepsilon}

\newcommand\ee{\nopagebreak\mbox{\ }\hfill$\diamond$}

\begin{document}

\title{A disintegration theorem for non-second-countable \'etale groupoids}

\date{July 7, 2026; revised July 16, 2026; minor changes August 7, 2026}

\author{Sergey Neshveyev}
\address{Department of Mathematics, University of Oslo, Norway}
\email{sergeyn@math.uio.no}

\begin{abstract}
We give a self-contained account of a version of Renault's disintegration theorem for (twisted) C$^*$-algebras of not necessarily Hausdorff \'etale groupoids that can be covered by countably many open bisections. As an application we discuss the $I$-norm, crossed product decompositions by inverse semigroup actions, and KMS-states and weights for general \'etale groupoid C$^*$-algebras.
\end{abstract}

\maketitle

\section*{Introduction}

Renault's disintegration theorem for second countable locally compact groupoids \cites{Rbook,MR0912813} is one of the foundational results in the theory of groupoid C$^*$-algebras. In full generality, its proof is quite long and requires some delicate function and measure theory. Various extensions, involving non-Hausdorff groupoids, twisted crossed products and Fell bundles~\cites{MR0912813,MR2547343,MR2446021,BHM}, make the picture even more intricate. As a result, the theorem has acquired the reputation of being such a heavy piece of machinery that it should be avoided whenever possible.

For \'etale groupoids, however, this reputation is largely undeserved. In this case most of the technical difficulties disappear, and the theorem becomes a relatively quick consequence of the disintegration theorem for representations of abelian C$^*$-algebras, which itself is little more than a version of the spectral theorem for bounded normal operators. Moreover, neither non-Hausdorffness nor twists introduce any substantial new difficulties, and the countability assumptions can be weakened considerably compared to the general case. None of these observations are of course revelatory. Nevertheless, both the literature and private discussions suggest that, even among experts, exactly what is involved in the disintegration theorem for \'etale groupoids and what the scope of its applicability is remain somewhat of a gray area. The purpose of this semi-expository note is to clarify these issues and present a few immediate applications.

\smallskip

We begin in Section~\ref{sec:integral} with a discussion of direct integrals and the disintegration of representations of abelian C$^*$-algebras. This material is classical and may be found in many textbooks; see, for example,~\cite{MR1451139}. However, since it is rarely presented in exactly the form needed here, we develop the theory from scratch in order to make it crystal clear what it involves. Fortunately, this requires relatively little space. The key point is that one has a satisfactory disintegration theory for representations of $C_0(X)$ for arbitrary Hausdorff locally compact spaces~$X$, provided one restricts attention to topologically countably generated Hilbert modules.

\smallskip

In Section~\ref{sec:disintegral} we prove a disintegration theorem for twisted \'etale groupoid C$^*$-algebras $C^*(\G,\LL)$, where the twist is given by a Fell line bundle~$\LL$. For pedagogical reasons we first treat the untwisted case in full detail, and then indicate the modifications needed for the twisted setting. In the untwisted case, our assumptions are that $\G$ can be covered by countably many open bisections and that the Hilbert modules under consideration are topologically countably generated. The twisted case requires slightly stronger hypotheses.

In recent years, \'etale groupoids that can be covered by countably many open bisections, but are not necessarily second countable, have attracted increasing attention. They arise naturally even when one starts from countable structures. For example, if $\Gamma$ is a nonamenable countable group, then its Furstenberg boundary $\partial_F\Gamma$ is not second countable~\cite{MR3652252}*{Corollary~3.17}. Consequently, the transformation groupoid $\Gamma\ltimes\partial_F\Gamma$ is not second countable, but it is clearly covered by countably many open bisections. As has also been observed by many, arbitrary \'etale groupoids can be expressed as increasing unions of open subgroupoids that admit countable covers by open bisections. Such groupoids therefore occupy a natural middle ground: the assumptions are broad enough to cover many important examples and to give information about arbitrary \'etale groupoids, while still often allowing the use of arguments that require countability.

\smallskip

In Section~\ref{sec:applications} we present several applications of the disintegration theorem. The first concerns Hahn's $I$-norm. This norm is used to define an $L^1$-type convolution algebra associated with a groupoid~$\G$, and $C^*(\G)$ is then defined as the C$^*$-envelope of this involutive Banach algebra~\cite{Rbook}. In the \'etale case, however, just as for discrete groups, one can instead define~$C^*(\G)$ simply as the universal C$^*$-completion of the convolution algebra $C_c(\G)$~\cite{MR2419901}. It is frustratingly nonobvious that these two constructions produce the same C$^*$-algebra, and although either definition suffices for most applications, the question of their equivalence has repeatedly arisen. In the second countable Hausdorff case this equivalence is a consequence of Renault's disintegration theorem. In full generality the equivalence has been recently established by Bardadyn, Kwa\'sniewski and McKee~\cite{MR4948043} using a disintegration theory for inverse semigroups of bisections that they developed. The version of Renault's theorem given here provides a shorter route to this result. Another related result that we prove is that for ample groupoids the C$^*$-algebra $C^*(\G)$ coincides with the universal C$^*$-completion of the Steinberg algebra of $\G$, which generalizes a result of Clark and Zimmerman~\cite{MR4750923}.

Our second application deals with crossed product decompositions of group\-oid C$^*$-algebras. The crossed product $C_0(X)\rtimes S$ for an action $S\curvearrowright X$ of an inverse semigroup was introduced by Sieben~\cite{MR1456588}. Shortly after, Paterson showed~\cite{MR1724106} that if $\G$ is a second countable \'etale groupoid and $S$ is a sufficiently rich countable inverse semigroup of open bisections of $\G$, then we have an isomorphism $C_0(\Gu)\rtimes S\cong C^*(\G)$. This result has been generalized in several ways \cites{MR4948043,MR2881538,MR3851326,MR3619758,MR2419901,MR1671944}, but it might still be difficult to extract from the literature the exact conditions under which the isomorphism holds. We show that, for arbitrary~$\G$ and~$S$, the isomorphism holds if and only if $\G$ decomposes as $S\ltimes\Gu$.

Our final application focuses on the structure of KMS-states and weights on $C^*(\G,\LL)$ for the dynamics defined by a $1$-cocycle $c\colon\G\to\R$. Renault showed~\cite{Rbook} that, when $\G$ is second countable and principal, there is a one-to-one correspondence between the KMS$_\beta$-states on $C^*(\G)$ and the quasi-invariant probability measures $\mu$ on $\Gu$ with Radon--Nikodym cocycle~$e^{-\beta c}$. When $\G$ is not principal, it was proved in~\cite{MR3138368} that the data describing a KMS$_\beta$-state consist, in addition to $\mu$, of an invariant Borel family of tracial states on the group C$^*$-algebras of the isotropy groups. This description has subsequently been generalized and refined in several directions~\cites{MR4222432,MR3861301,MR4592883,MR4490951,MR4411323}. Since all of these works rely on Renault's disintegration theorem, second countability has remained a standing assumption throughout. We show that it can be replaced by the weaker hypotheses used in Section~\ref{sec:disintegral}. Along the way, and as an application, we also prove some results for general~$\G$.

\medskip\noindent
{\bf Acknowledgement.} It is a pleasure to thank Bartosz Kwa\'sniewski for an inspiring correspondence and his comments on this work.

\bigskip

\section{Direct integral decompositions}\label{sec:integral}

\subsection{Measurable fields of Hilbert spaces} \label{ssec:measure}
Assume $(X,\BB)$ is a measurable space.

\begin{defn}\label{def:measurable-fields}
A \emph{measurable field of separable Hilbert spaces} over $(X,\BB)$ is a collection $\HH=(H_x)_{x\in X}$ of such spaces together with a subspace $\Gamma\subset\prod_{x\in X}H_x$ of sections of $\HH$ satisfying the following properties:
\begin{enumerate}
  \item[(i)] for every $\xi=(\xi_x)_{x\in X}\in\Gamma$, the function $x\mapsto\|\xi_x\|$ is measurable;
  \item[(ii)] if a section $\zeta=(\zeta_x)_{x\in X}$ of $\HH$ has the property that the function $x\mapsto (\zeta_x,\xi_x)$ is measurable for all $\xi\in\Gamma$, then $\zeta\in\Gamma$;
  \item[(iii)] there is a countable subset of $\Gamma$ such that its image in $H_x$ is dense for all $x\in X$.
\end{enumerate}

Two measurable fields $(\HH,\Gamma)$ and $(\HH',\Gamma')$ are called \emph{unitarily isomorphic} if there exists a collection $(U_x)_{x\in X}$ of unitary operators $U_x\colon H_x\to H_x'$ such that a section $\xi$ of $\HH$ lies in $\Gamma$ if and only if the section $(U_x\xi_x)_x$ lies in $\Gamma'$.
\end{defn}

We will denote a measurable field by just one letter $\HH$, instead of a pair $(\HH,\Gamma)$, and call the elements of $\Gamma$ \emph{measurable sections} of $\HH$. The adjective \emph{separable} in the name is usually omitted, but since the precise role of various countability conditions is important for us, we prefer to keep~it. Another name for $\HH$ used in the literature is a \emph{measurable Hilbert bundle}.

\begin{example}
Take a separable Hilbert space $H$ and define a measurable field $\HH=(H_x)_{x\in X}$ by letting $H_x:=H$ and declaring a section $(\xi_x)_x$ to be measurable iff the functions $x\mapsto (\xi_x,\zeta)$ are measurable for all $\zeta\in H$. Condition (i) in Definition~\ref{def:measurable-fields} is checked by choosing an orthonormal basis $(e_k)_k$ in $H$ and using the Parseval identity, the other two conditions are obvious. Note also that a section $(\xi_x)_x$ is measurable if and only if the Fourier coefficients $x\mapsto (\xi_x,e_k)$ are measurable. We call~$\HH$ the \emph{constant field with fiber} $H$ or a \emph{constant field of rank} $\dim H$.

Let us also mention that it is not difficult to show that the norm and weak topologies on $H$ define the same Borel $\sigma$-algebra $\BB_H$, and then a section $(\xi_x)_x$ is measurable if and only if the map $(X,\BB)\to (H,\BB_H)$, $x\mapsto\xi_x$, is measurable. \ee
\end{example}

Let us discuss a few simple consequences of the definition and introduce some notation. Given a measurable field  $\HH=(H_x)_{x\in X}$  of separable Hilbert spaces and two measurable sections~$\xi$ and~$\zeta$, the polarization identity implies that the function $x\mapsto(\xi_x,\zeta_x)$ is measurable.

If we are given a countable collection of  measurable fields $\HH^{(n)}=(H_x^{(n)})_{x\in X}$  of separable Hilbert spaces, then we can define in the obvious way their direct sum $\bigoplus_n\HH^{(n)}$ such that its section is measurable if and only if its $\HH^{(n)}$-component is measurable for each~$n$.

If $\xi$ is a measurable section of a measurable field  $\HH=(H_x)_{x\in X}$  of separable Hilbert spaces and $f\colon X\to\C$ is a measurable function, then by condition (ii) in Definition~\ref{def:measurable-fields} the section $f\xi=(f(x)\xi_x)_x$ is measurable as well. Given a measurable subset $Y\subset X$, we can consider the field $(H_y)_{y\in Y}$ and declare its section $\xi$ to be measurable iff it is the restriction of a measurable section of $\HH=(H_x)_{x\in X}$, or equivalently, iff the section $\un_Y\xi$ of $\HH$ is measurable. We denote the measurable field we thus get by~$\HH|_Y$. If $X=\bigsqcup_nY_n$ is a countable partition of $X$ into measurable sets, then a section~$\xi$ of~$\HH$ is measurable if and only if $\xi|_{Y_n}$ is a measurable section of $\HH|_{Y_n}$ for all~$n$. Together with the following result this often allows one to reduce general arguments to the case of constant fields.

\begin{prop}\label{prop:trivialization}
Assume $\HH=(H_x)_x$ is a measurable field of separable Hilbert spaces over a measurable space $(X,\BB)$. For every $n\in\Z_+\cup\{\infty\}$, denote by $X_n\subset X$ the subset of points $x$ such that $\dim H_x=n$. Then, for every $n$, the set $X_n$ is measurable and $\HH|_{X_n}$ is unitarily isomorphic to a constant field of rank $n$.
\end{prop}

\bp
Take a sequence of measurable sections $\zeta^{(k)}$ ($k\in\N$) such that its image in $H_x$ is dense for all $x\in X$. The result is obtained by applying fiber-wise the Gram--Schmidt orthogonalization procedure to this sequence. In detail, denote by $P^{(k)}_x$ the orthogonal projection $H_x\to\operatorname{span}\{\zeta^{(1)}_x,\dots,\zeta^{(k)}_x\}$. Let also $P^{(0)}_x:=0$. For $k\ge1$, we define $\xi^{(k)}_x:=0$ if $P^{(k-1)}_x\zeta^{(k)}_x=\zeta^{(k)}_x$, otherwise we let
$$
\xi^{(k)}_x:=\frac{\zeta^{(k)}_x-P^{(k-1)}_x\zeta^{(k)}_x}{\|\zeta^{(k)}_x-P^{(k-1)}_x\zeta^{(k)}_x\|}.
$$
For every $x$, the following properties hold: the vectors $\xi^{(k)}_x$ are mutually orthogonal, they span a dense subspace of $H_x$, and, for every $k\in\N$, we have either $\xi^{(k)}_x=0$ or $\|\xi^{(k)}_x\|=1$. Since $P^{(k)}_x=\sum^{k}_{l=1}(\cdot,\xi^{(l)}_x)\xi^{(l)}_x$, by induction on $k$ we see that the sections $\xi^{(k)}$ are measurable.
This already implies that the function $x\mapsto\dim H_x=\sum^\infty_{k=1}\|\xi^{(k)}_x\|$ is measurable, hence the sets $X_n$ are measurable.

Let us now fix $n\in\N\cup\{\infty\}$ and a Hilbert space $H$ with an orthonormal basis $(e_k)^n_{k=1}$. Let $N_k(x):=\operatorname{rank} P^{(k)}_x=\sum^k_{l=1}\|\xi^{(l)}_x\|$. For every $x\in X_n$, define a unitary operator $U_x\colon H_x\to H$ by omitting the zero vectors in the sequence $(\xi^{(k)}_x)^\infty_{k=1}$ and mapping the orthonormal basis of~$H_x$ that we get onto $(e_k)^n_{k=1}$. In other words, $U_x\xi^{(k)}_x=e_{N_k(x)}$ if $\xi^{(k)}_x\ne0$. Since the functions~$N_k$ are measurable and $\xi_x=\sum^\infty_{k=1}(\cdot,\xi^{(k)}_x)\xi^{(k)}_x$ for any section $\xi$ of $\HH|_{X_n}$, the family $(U_x)_{x\in X_n}$ maps every measurable section of $\HH|_{X_n}$ into a measurable section of the constant field with fiber $H$ over $X_n$.

The inverse operator $U_x^{-1}$ maps $e_k$ into $\xi^{(l_k(x))}_{x}$, where $l_k(x):=\min\{l\ge k: N_l(x)=k\}$. Since the functions $l_k$ are measurable, we similarly conclude that the family $(U_x^{-1})_{x\in X_n}$ maps every measurable section of the constant field with fiber $H$ into a measurable section of $\HH|_{X_n}$. Thus, $(U_x)_{x\in X_n}$ defines the required isomorphism.
\ep

Therefore a measurable field of separable Hilbert spaces contains no more information than a measurable function $X\to \Z_+\cup\{\infty\}$. The point of this notion, however, is that it provides a convenient framework for analyzing measurable families of operators.

Assume we are given two measurable fields of separable Hilbert spaces $\HH=(H_x)_x$ and $\HH'=(H'_x)_x$ over $(X,\BB)$ and, for every $x\in X$, a bounded linear operator $T_x\colon H_x\to H'_x$. We say that $(T_x)_x$ is a \emph{measurable family of operators} if the section $(T_x\xi_x)_x$ of $\HH'$ is measurable for every measurable section $\xi=(\xi_x)_x$ of $\HH$. We say that the family is bounded if $\sup_x\|T_x\|<\infty$.

It is clear that if $(T_x\colon H_x\to H'_x)_x$ and $(S_x\colon H'_x\to H''_x)_x$ are two measurable families of operators, then the family $(S_xT_x)_x$ is measurable as well. Also, by condition (ii) in Definition~\ref{def:measurable-fields}, if $(T_x)_x$ is measurable, then $(T^*_x)_x$ is measurable as well. In particular, any measurable family of unitary operators defines a unitary isomorphism in the sense of Definition~\ref{def:measurable-fields}. Therefore the last paragraph of the proof of Proposition~\ref{prop:trivialization} was in fact redundant.

When $\HH$ and $\HH'$ are constant fields with fibers $H$ and $H'$, resp., a measurable family of operators is a map $X\to B(H,H')$, $x\mapsto T_x$, that is weakly measurable, meaning that the functions $x\mapsto(T_x\xi,\zeta)$ are measurable for all $\xi\in H$ and $\zeta\in H'$.

\subsection{Direct integrals and decomposable operators}

Assume $\HH=(H_x)_x$ is a measurable field of separable Hilbert spaces over a measurable space $(X,\BB)$. Assume in addition that $\mu$ is a $\sigma$-finite measure on $(X,\BB)$. A measurable section $\xi$ of $\HH$ is called \emph{square-integrable} if
$$
\int_X\|\xi_x\|^2d\mu(x)<\infty.
$$
Identifying square-integrable sections that coincide $\mu$-a.e., we get a Hilbert space $\int^\oplus_XH_xd\mu(x)$ of such sections with inner product
$$
(\xi,\zeta):=\int_X(\xi_x,\zeta_x)d\mu(x).
$$
It is called the \emph{direct integral} of $(H_x)_x$ over $(X,\BB,\mu)$. The completeness of $\int^\oplus_XH_xd\mu(x)$ can be checked in the same way as for the $L^2$-spaces. Alternatively, by taking a direct sum decomposition of $\HH$ and by partitioning $X$ one can reduce the problem to the case of a constant field $\HH$ of rank one, in which case $\int^\oplus_XH_xd\mu(x)$ becomes unitarily isomorphic to $L^2(X,d\mu)$.

Assume we are given another measurable field $\HH'=(H'_x)_x$ of separable Hilbert spaces and a bounded measurable family of operators $T_x\colon H_x\to H_x'$. We then get a well-defined operator
$$
T=\int^\oplus_XT_xd\mu(x)\colon \int^\oplus_XH_xd\mu(x)\to \int^\oplus_XH'_xd\mu(x)
$$
mapping a square-integrable section $\xi$ of $\HH$ into $(T_x\xi_x)_x$. It is clear that this operator is bounded, with $\|T\|\le\sup_x\|T_x\|$; in fact, as we will see shortly, $\|T\|=\operatorname{ess.supp}_x\|T_x\|$. Such operators $\int^\oplus_XT_xd\mu(x)$ are called \emph{decomposable}.

We have canonical representations of $L^\infty(X,\mu)$ on $\int^\oplus_XH_xd\mu(x)$ and $\int^\oplus_XH'_xd\mu(x)$ by multiplication operators, also called \emph{diagonal operators}. By slightly abusing notation, we denote these representations by one letter $m$. Thus, $m(f)\xi=f\xi$ for any square-integrable section $\xi$ of~$\HH$ or~$\HH'$. Any decomposable operator $T\colon \int^\oplus_XH_xd\mu(x)\to \int^\oplus_XH'_xd\mu(x)$ satisfies $Tm(f)=m(f)T$ for all $f\in L^\infty(X,\mu)$. As the following result shows, this is a characteristic property of decomposable operators.

\begin{prop}\label{prop:decomposable-operators}
Assume $(X,\BB,\mu)$ is a $\sigma$-finite measure space, $\HH=(H_x)_x$ and $\HH'=(H'_x)_x$ are measurable fields of separable Hilbert spaces over $(X,\BB)$, and $T\colon\int^\oplus_XH_xd\mu(x)\to\int^\oplus_X H'_xd\mu(x)$ is a bounded operator such that $Tm(f)=m(f)T$ for all $f\in L^\infty(X,\mu)$. Then there exists a bounded measurable family of operators $T_x\colon H_x\to H'_x$ such that $T=\int^\oplus_XT_xd\mu(x)$, and if $(T'_x)_x$ is another family with the same property, then $T_x=T'_x$ for $\mu$-a.e.~$x$.
\end{prop}

\bp
Since we can view $T$ as an operator on $\int^\oplus_XH_xd\mu(x)\bigoplus \int^\oplus_XH'_xd\mu(x)=\int^\oplus_X(H_x\oplus H'_x)d\mu(x)$, we may assume that $\HH=\HH'$. If $A\in\BB$, then the space $m(\un_A)\int^\oplus_XH_xd\mu(x)=\int^\oplus_AH_xd\mu(x)$ is $T$-invariant. In view of Proposition~\ref{prop:trivialization} and $\sigma$-finiteness of $\mu$, it follows that without loss of generality we may also assume that $\mu$ is finite and $\HH$ is a constant field with fiber $H$.

Fix an orthonormal basis $(e_k)_k$ in $H$. Viewing $e_k$ as a section of $\HH$ and fixing representatives of elements of $\int^\oplus_XH_xd\mu(x)$,
we get measurable sections $x\mapsto (Te_k)_x\in H$ of $\HH$. Consider the measurable functions $x\mapsto T_x(k,l)$ defined by $T_x(k,l):=((Te_l)_x,e_k)$. We thus get matrices $T_x:=(T_x(k,l))_{k,l}$. We claim that $\|T_x\|\le\|T\|$ for $\mu$-a.e.~$x$.

In order to prove this, choose a dense sequence $(s_i)^\infty_{i=1}$ of vectors in the intersection of the linear span of the basis vectors $e_k$ with the unit ball of $H$. Then the vectors $T_xs_i$ are well-defined and $\|T_x\|=\sup_i\|T_xs_i\|$, so it suffices to show that the measurable functions $h_i$ defined by $h_i(x):=\|T_xs_i\|^2$ satisfy $\|h_i\|_\infty\le\|T\|^2$. For $A\in\BB$, the section $\un_A s_i$ is square-integrable and $T(\un_A s_i)=Tm(\un_A)s_i=m(\un_A)Ts_i$, so that $(T(\un_As_i))_x=\un_A(x)T_xs_i$ a.e. Hence
$$
\int_A h_id\mu=\|T(\un_As_i)\|^2\le\|T\|^2\|\un_A s_i\|^2\le\|T\|^2\mu(A).
$$
This implies that $\|h_i\|_\infty\le\|T\|^2$, proving our claim.

By modifying the sections $x\mapsto(Te_k)_x$ on a set of measure zero, we may assume that $\|T_x\|\le\|T\|$ for all $x$. Since the functions $x\mapsto T_x(k,l)$ are measurable, the family $(T_x)_x$ is measurable. Consider the operator $\tilde T:=\int^\oplus_XT_xd\mu(x)$. By construction we have $Te_k=\tilde Te_k$ for all $k$, and hence $T(fe_k)=m(f)Te_k=\tilde T (fe_k)$ for all $f\in L^\infty(X,\mu)\subset L^2(X,d\mu)$. Since such vectors $fe_k$ span a dense subspace of $\int^\oplus_XH_xd\mu(x)$, it follows that $T=\tilde T$.

It remains to show that if $T=\int^\oplus_X T'_xd\mu(x)$ is another decomposition, then $T_x=T'_x$ for $\mu$-a.e.~$x$. In fact, a stronger property holds: for any decomposable operator $T=\int^\oplus_X T_xd\mu(x)$, we have $\|T\|=\operatorname{ess.supp}_x\|T_x\|$. The inequality $\le$ is obvious. The opposite inequality holds by the claim in the first part of the proof.
\ep

Since we obviously have $(\int^\oplus_XT_xd\mu(x))^*=\int^\oplus_XT_x^*d\mu(x)$, the essential uniqueness of a decomposition $T=\int^\oplus_XT_xd\mu(x)$ gives the following result.

\begin{cor}\label{cor:unitary-selfadjoint}
A decomposable operator $T=\int^\oplus_XT_xd\mu(x)$ is self-adjoint if and only if $T_x$ is self-adjoint for $\mu$-a.e.~$x$, and it is unitary if and only if $T_x$ is unitary for $\mu$-a.e.~$x$.
\end{cor}

\subsection{Representations of abelian C\texorpdfstring{$^*$}{*}-algebras} \label{ssec:abelian}

Assume now that $X$ is a Hausdorff locally compact space. Denote by $\BB_X$ the Borel $\sigma$-algebra on $X$. By a \emph{Borel field of separable Hilbert spaces} over $X$ we mean a measurable field $\HH$ of such spaces over $(X,\BB_X)$. Correspondingly, we then talk about Borel sections and families of operators instead of measurable ones.

If $\mu$ is a $\sigma$-finite regular Borel measure on $X$, we can consider the representation~$\pi$ of~$C_0(X)$ by diagonal operators on $\int^\oplus_X H_xd\mu(x)$. If $\HH$ is a constant field of rank one, then this representation is unitarily equivalent to the standard representation by multiplication operators on $L^2(X,d\mu)$, and then by $\sigma$-finiteness and regularity of $\mu$ we conclude that the representation~$\pi$ is cyclic. Proposition~\ref{prop:trivialization} implies that for general~$\HH$ the representation $\pi$ is a direct sum of countably many cyclic representations. Equivalently, the $C_0(X)$-module $\int^\oplus_X H_xd\mu(x)$ is topologically countably generated. The following result shows that all topologically countably generated Hilbert modules arise this way.

\begin{prop}\label{prop:abelian}
Assume $X$ is a Hausdorff locally compact space, $\pi\colon C_0(X)\to B(H)$ is a nondegenerate representation such that the $C_0(X)$-module $H$ is topologically countably generated. Then there exists a regular Borel probability measure $\mu$ on $X$ and a Borel field $\HH=(H_x)_x$ of separable Hilbert spaces over $X$ such that $H_x\ne0$ for $\mu$-a.e.~$x$ and $\pi$ is unitarily equivalent to the representation of~$C_0(X)$ by diagonal operators on $\int^\oplus_X H_xd\mu(x)$. If $(\mu',\HH')$ is another pair with the same properties, then $\mu\sim\mu'$ and the fields~$\HH$ and~$\HH'$ are unitarily isomorphic outside a Borel subset of $X$ of $\mu$-measure zero.
\end{prop}

\bp
Recall first that any unit vector $\xi\in H$ defines a state $(\pi(\cdot)\xi,\xi)$ on $C_0(X)$ and hence a unique regular Borel probability measure $\nu_\xi$ on $X$ such that
$$
(\pi(f)\xi,\xi)=\int_X f(x)d\nu_\xi(x)\quad\text{for all}\quad f\in C_0(X).
$$

Now, by assumption, we can decompose $\pi$ into a direct sum of countably many cyclic representations $\pi_n\colon C_0(X)\to B(H^{(n)})$. For every $n$, fix a unit cyclic vector $\xi_n\in H^{(n)}$ and consider the corresponding probability measure $\mu_n:=\nu_{\xi_n}$. Choose scalars $\lambda_n>0$ such that $\sum_n\lambda_n=1$ and put $\mu:=\sum_n\lambda_n\mu_n$. In other words, $\mu=\nu_\xi$ for $\xi:=\sum_n\lambda_n^{1/2}\xi_n$. Each representation $\pi_n$ is equivalent to the representation of $C_0(X)$ by multiplication operators on $L^2(X,d\mu_n)$. Consider the Borel field $\HH^{(n)}=(H^{(n)}_x)_x$ such that $H^{(n)}_x=\C$ for  $x\in A_n:=\{y:\frac{d\mu_n}{d\mu}(y)>0\}$ and $H^{(n)}_x=0$ for $x\in X\setminus A_n$, so that $\HH^{(n)}|_{A_n}$ and $\HH^{(n)}|_{X\setminus A_n}$ are constant fields of rank $1$ and $0$, resp. Then the map $f\mapsto \big(\frac{d\mu_n}{d\mu}\big)^{1/2}f$ defines a unitary isomorphism $L^2(X,d\mu_n)\cong\int^\oplus_XH^{(n)}_xd\mu(x)$ intertwining the actions of $C_0(X)$ by multiplication operators. It follows that $\HH:=\bigoplus_n\HH^{(n)}$ is the required Borel field.

Before we turn to the uniqueness part, let us observe that if $(\mu,\HH)$ is as in the statement of the proposition and $\pi$ is identified with the representation of $C_0(X)$ on $\int^\oplus_X H_xd\mu(x)$, then the following properties hold:
\begin{itemize}
  \item[(a)] we have $\nu_\xi\ll\mu$ for all unit vectors $\xi\in H$;
  \item[(b)] $\mu=\nu_\xi$ for some unit vector $\xi\in H$;
  \item[(c)] the representation $m$ of $L^\infty(X,\mu)$ by diagonal operators on $H=\int^\oplus_X H_xd\mu(x)$ is the unique normal extension of $\pi$.
\end{itemize}
Property (a) holds, because for any square-integrable section $\xi=(\xi_x)_x$ we have $\int_X f\,d\nu_\xi=\int_X f(x)\|\xi_x\|^2d\mu$. Property (b) follows by choosing a Borel section $\xi=(\xi_x)_x$ such that $\|\xi_x\|=1$ for $\mu$-a.e.~$x$, which is possible by Proposition~\ref{prop:trivialization} and the assumption that $H_x\ne0$ a.e. Finally, since any representation is determined by its matrix coefficients, property~(c) is equivalent to the following two statements: (c$'$) the linear functionals $(m(\cdot)\xi,\zeta)$ are normal, that is, they lie in $L^1(X,d\mu)\subset L^\infty(X,\mu)^*$, for all square-integrable sections~$\xi$ and~$\zeta$, and (c$''$) the image of $C_0(X)$ in $L^\infty(X,\mu)$ is weakly$^*$ dense. The first statement is true, since $(m(f)\xi,\zeta)=\int_X f(x)(\xi_x,\zeta_x)d\mu(x)$, while the second statement follows from regularity of~$\mu$.

Now, assume we have another pair $(\mu',\HH')$ as in the statement of the proposition. Properties~(a) and~(b) applied to $\mu$ and $\mu'$ imply that these measures must be equivalent. The map $\xi\mapsto \big(\frac{d\mu'}{d\mu}\big)^{1/2}\xi$ defines a unitary isomorphism $\int^\oplus_X H'_xd\mu'(x)\to \int^\oplus_X H'_xd\mu(x)$ intertwining the actions of $C_0(X)$ by diagonal operators. Therefore the pair $(\mu,\HH')$ also defines a decomposition of~$\pi$ as in the statement of the proposition. Then the identity operator on $H$ can be viewed as a unitary operator $U\colon \int^\oplus_X H_xd\mu(x)\to \int^\oplus_X H'_xd\mu(x)$ intertwining the actions of $C_0(X)$ by diagonal operators. By property (c) above, this operator must also intertwine the actions of $L^\infty(X,\mu)$. By Proposition~\ref{prop:decomposable-operators} it follows that $U$ is decomposable, so $U=\int^\oplus_XU_xd\mu(x)$. By Corollary~\ref{cor:unitary-selfadjoint} the operators $U_x\colon H_x\to H_x'$ are unitary for $\mu$-a.e.~$x$, so they define the required isomorphism of $\HH$ onto $\HH'$ outside a set of measure zero.
\ep

The essential uniqueness statement in Proposition~\ref{prop:decomposable-operators} implies that the last part of the above proof gives the following more precise result.

\begin{cor}\label{cor:abelian}
Assume $X$ is a Hausdorff locally compact space, $\mu$ and $\mu'$ are regular Borel probability measures on $X$, and $\HH=(H_x)_x$ and $\HH'=(H'_x)_x$ are Borel fields of separable Hilbert spaces over $X$ such that $H_x\ne0$ for $\mu$-a.e.~$x$ and $H'_x\ne0$ for $\mu'$-a.e.~$x$. Assume $U\colon\int^\oplus_X H_xd\mu(x)\to \int^\oplus_X H'_xd\mu'(x)$ is a unitary operator intertwining the representations of $C_0(X)$ by diagonal operators. Then $\mu\sim\mu'$ and there is a Borel family of operators $U_x\colon H_x\to H'_x$ such that $U_x$ is unitary for $\mu$-a.e.~$x$ and $U\xi=\big(\frac{d\mu}{d\mu'}(x)^{1/2}U_x\xi_x\big)_x$ for all $\xi\in\int^\oplus_X H_xd\mu(x)$. If $(U'_x)_x$ is another Borel family with the same properties, then $U_x=U'_x$ for $\mu$-a.e.~$x$.
\end{cor}

\bigskip

\section{Disintegration theorems}\label{sec:disintegral}

\subsection{Representations of groupoid C\texorpdfstring{$^*$}{*}-algebras}

Let $\G$ be a locally compact \'etale groupoid with unit space $\Gu$, the range map $r\colon g\mapsto gg^{-1}$ and the source map $s\colon g\mapsto  g^{-1}g$. We assume that the unit space $\Gu$ is Hausdorff in the relative topology, but we do not require the entire groupoid $\G$ to be Hausdorff. \'Etaleness means by definition that the maps $r,s\colon\G\to\Gu$ are local homeomorphisms, which then implies that~$\Gu$ is open in $\G$. For a subset $A\subset\Gu$, we write $\G^A$ for $r^{-1}(A)$ and $\G_A$ for $s^{-1}(A)$. The isotropy groups are $\Gxx:=\G^x\cap\G_x$.

For every open Hausdorff subset $U\subset\G$, consider the space $C_c(U)$ of continuous compactly supported functions on $U$. We view the elements of $C_c(U)$ as functions on $\G$ by extending them by zero. Note that if $\G$ is non-Hausdorff, then these extensions are not necessarily continuous. We then define $C_c(\G)$ as the linear span of all such spaces $C_c(U)$; other notation for this space one can find in the literature are $C_c(\G)_{0}$,  $\mathscr{C}_c(\G)$, $\mathcal{C}(\G)$. A partition of unity argument shows that if $(U_i)_i$ is an open cover of $\G$ consisting of Hausdorff sets, then $C_c(\G)$ is spanned by the spaces~$C_c(U_i)$.

Recall that a subset $W\subset\G$ is called a bisection if the maps $r|_W$ and $s|_W$ are injective. If $W\subset\G$ is an open bisection, then the maps $r\colon W\to r(W)$ and $s\colon W\to s(W)$ are homeomorphisms. In particular, every such bisection $W$ is Hausdorff. Since $\G$ is covered by open bisections, the space $C_c(\G)$ is spanned by the spaces $C_c(W)$.

The space $C_c(\G)$ is a $*$-algebra with convolution product
\begin{equation*} \label{eprod}
(f_{1}*f_{2})(g) := \sum_{h \in \G^{r(g)}} f_{1}(h) f_{2}(h^{-1}g)= \sum_{h \in \G_{s(g)}} f_{1}(gh^{-1}) f_{2}(h)
\end{equation*}
and involution by $f^{*}(g):=\overline{f(g^{-1})}$.

\begin{lemma}\label{lem:unibound}
For every $f\in C_c(\G)$, there is a constant $C_f\ge0$ depending only on $f$ such that if~$\pi$ is a representation of $C_c(\G)$ on a pre-Hilbert space $H$, then $\|\pi(f)\|\le C_f$. More precisely, if $W\subset\G$ is an open bisection and $f\in C_c(W)$, then
\begin{equation}\label{eq:upper-bound}
\|\pi(f)\|\le\|f\|_\infty.
\end{equation}
\end{lemma}

Here by a representation we mean an assignment to every $f\in C_c(\G)$ of a linear operator $\pi(f)\colon H\to H$, which is not a priori assumed to be bounded, such that $\pi$ is an algebra homomorphism and $(\pi(f)\xi,\zeta)=(\xi,\pi(f^*)\zeta)$ for all $f\in C_c(\G)$ and $\xi,\zeta\in H$.

\bp
Since $C_c(\G)$ is spanned by the spaces $C_c(W)$ for open bisections $W$, it suffices to prove the last statement of the lemma.

Let us check first that if $h\in C_c(\Gu)$ is such that $0\le h\le 1$, then $\|\pi(h)\|\le1$. Since $(h-h^2)^{1/2}\in C_c(\Gu)$, for any $\xi\in H$ we have $((\pi(h)-\pi(h)^2)\xi,\xi)\ge0$. Then
$$
\|\pi(h)\xi\|^2\le(\pi(h)\xi,\xi)\le\|\pi(h)\xi\|\,\|\xi\|,
$$
and hence $\|\pi(h)\xi\|\le\|\xi\|$. Thus, $\|\pi(h)\|\le1$.

Now, take an open bisection $W$ and a nonzero $f\in C_c(W)$. Then $h:=\|f\|_\infty^{-2}f^**f$ satisfies $h\in C_c(\Gu)$ and $0\le h\le1$, hence $\|\pi(h)\|\le1$ and therefore $\|\pi(f)\|^2=\|\pi(f^**f)\|\le\|f\|_\infty^2$.
\ep

This lemma implies that we get a C$^*$-seminorm on $C_c(\G)$ defined by
$$
\|f\|:=\sup_\pi\|\pi(f)\|,
$$
where the supremum is taken over all representations $\pi$ of $C_c(\G)$ by bounded operators on Hilbert spaces, or equally well, over all representations by arbitrary operators on pre-Hilbert spaces. This is in fact a norm dominating the supremum-norm, since for every point $x\in \Gu $ we have a representation $\rho_{x}\colon C_{c}(\G) \to B(\ell^{2}(\G_{x}))$ defined~by
\begin{equation}\label{eq:rhox2}
(\rho_x(f)\xi)(g) :=
\sum_{h\in \G^{r(g)}}f(h) \xi(h^{-1}g),
\end{equation}
and then for all $g\in\G_x$ we have $(\rho_x(f)\delta_x,\delta_g)=f(g)$. The full C$^*$-algebra $C^*(\G)$ of $\G$ is defined as the completion of $C_c(\G)$ in this norm.

Since the C$^*$-norm on $C_c(\G)$ dominates the supremum-norm and at the same time is dominated by it on $C_c(W)$ for any open bisection $W$ by~\eqref{eq:upper-bound}, we conclude that
\begin{equation}\label{eq:norm-equality}
\|f\|=\|f\|_\infty
\end{equation}
for all $f\in C_c(W)$. In particular, taking $W=\Gu$ we see that the closure of~$C_c(\Gu)$ in~$C^*(\G)$ is naturally identified with~$C_0(\Gu)$.

\smallskip

Before turning to disintegration of representations of $C^*(\G)$, let us discuss a few more aspects of the theory.

The open bisections of $\G$ form an inverse semigroup $\BisG$ under multiplication. If $S\subset\BisG$ is a sufficiently rich inverse subsemigroup, then $C^*(\G)$ can be identified with the crossed product $C_0(\Gu)\rtimes S$. We will discuss this in detail in Section~\ref{ssec:inverse}. For now, it will suffice to understand how representations of $C^*(\G)$ give rise to representations of~$\BisG$.

By a representation of an inverse semigroup $S$ on a Hilbert space~$H$ we mean an assignment to every $s\in S$ of a partial isometry $u_s\in B(H)$ such that $u_{st}=u_su_t$ and $u_s^*=u_{s^{-1}}$. 

\begin{lemma}[{cf.~\cite{MR1456588}*{Proposition~5.6}}]\label{lem:bis-rep}
Assume we are given a representation $\pi\colon C^*(\G)\to B(H)$. For every open set $U\subset\Gu$, denote by $p_U$ the projection onto $\overline{\pi(C_0(U))H}$. Then, for every $W\in\BisG$, there is a unique partial isometry $u_W\in B(H)$ with initial projection~$p_{s(W)}$ such that for all $f\in C_c(s(W))$ we have $u_W\pi(f)=\pi(\un_W*f)$. Moreover, $u_W\in\pi(C^*(\G))''$, the range projection of $u_W$ is $p_{r(W)}$, and the map $W\mapsto u_W$ is a representation of the inverse semigroup~$\BisG$.
\end{lemma}

Note that although $\un_W\notin C_c(\G)$ in general, the convolution $\un_W*f$ is a well-defined element of $C_c(W)\subset C_c(\G)$ for $f\in C_c(s(W))$, as $(\un_W*f)(g)=f(s(g))$ for $g\in W$.

\bp
The uniqueness is clear. For every open set $U\subset\Gu$, the net $\Lambda_U:=\{h\in C_c(U):0\le h\le1\}$ is an approximate unit in $C_0(U)$. As is well-known, then $\pi(h)\to_{h\in\Lambda_U} p_U$ in the strong operator topology. For similar reasons, for every $W\in\BisG$, the net $(\pi(\un_W*h))_{h\in\Lambda_{s(W)}}$ converges strongly to a partial isometry $u_W$ with initial projection $p_{s(W)}$, since if $f\in C_c(s(W))$ and $\xi\in H$, then $\pi(\un_W*h)\pi(f)\xi=\pi(\un_W*(hf))\xi=\pi(\un_W*f)\xi$ for all $h\in\Lambda_{s(W)}$ large enough, and if $f_1,f_2\in C_c(s(W))$ and $\xi_1,\xi_2\in H$, then
$$
(\pi(\un_W*f_1)\xi_1,\pi(\un_W*f_2)\xi_2)=(\pi(\bar f_2*\un_{W^{-1}}*\un_W*f_1)\xi_1,\xi_2)=(\pi(f_1)\xi_1,\pi(f_2)\xi_2).
$$

We then have $\pi(\un_W*h)^*\to_{h\in \Lambda_{s(W)}} u_W^*$ in the weak operator topology. Since
$$
(\un_W*h)^*=\un_{W^{-1}}*(\un_{W}*h*\un_{W^{-1}})
$$
and $\{\un_{W}*h*\un_{W^{-1}}:h\in\Lambda_{s(W)}\}=\Lambda_{r(W)}=\Lambda_{s(W^{-1})}$, it follows that $u_W^*=u_{W^{-1}}$. In particular, the range projection of $u_W$ is~$p_{r(W)}$. 

It remains to show that the map $W\mapsto u_W$ is multiplicative. This follows by observing that if $V,W\in\BisG$, $h_1\in\Lambda_{s(V)}$ and $h_2\in\Lambda_{s(W)}$, then
$$
\un_V*h_1*\un_W*h_2=\un_{VW}*\big((\un_{W^{-1}}*h_1*\un_W)h_2\big),
$$
the functions $(\un_{W^{-1}}*h_1*\un_W)h_2$ lie in $C_c(s(VW))$ (equivalently, the functions $h_1(\un_{W}*h_2*\un_{W^{-1}})$ lie in $C_c(s(V)\cap r(W))$) and form an approximate unit in $C_0(s(VW))$, so that the net $(\pi(\un_V*h_1*\un_W*h_2))_{h_1,h_2}$ converges strongly to both $u_Vu_W$ and $u_{VW}$.
\ep

Every open bisection $W$ defines a homeomorphism $T_W\colon s(W)\to r(W)$, the composition of $r$ with the inverse of $s\colon W\to s(W)$. A $\sigma$-finite regular Borel measure $\mu$ on $\Gu$ is called $\G$-\emph{quasi-invariant} if the pushforward of $\mu|_{s(W)}$ under $T_W$ is equivalent to $\mu|_{r(W)}$ for all $W\in\BisG$. A Borel function $D\colon \G\to(0,+\infty)$ is then called the \emph{Radon--Nikodym cocycle} of $\mu$ if $\frac{d\left((T_{W})_*\mu\right)}{d\mu}(x)=D(g^x_W)^{-1}$ for $\mu$-a.e.~$x\in r(W)$ and all Borel bisections $W$, where $g^x_W$ is the unique element of $W\cap\G^x$. As will follow from our considerations in Section~\ref{ssec:KMS}, in general we can only guarantee that $D$ can be chosen to be locally Borel, but in this section we are primarily interested in groupoids that can be covered by countably many open bisections, in which case existence of a Borel Radon--Nikodym cocycle is straightforward. Moreover, for such groupoids we get that
\begin{equation}\label{eq:RN}
\int_{\Gu}\sum_{g\in\G_x}f(g)d\mu(x)=\int_{\Gu}\sum_{g\in\G^x}f(g)D(g)^{-1}d\mu(x)
\end{equation}
for all positive Borel functions $f$ on $\G$. Indeed, if $f=\un_A$ for a Borel subset $A$ of an open bisection~$W$, then this identity holds by definition, and then for general $f$ the identity is obtained by an application of the monotone convergence theorem.

\begin{defn}
Assume $\HH=(H_x)_{x\in\Gu}$ is a Borel field of separable Hilbert spaces over $\Gu$.
By a \emph{Borel action} of $\G$ on $\HH$ we mean an assignment to every $g\in\G$ of a unitary operator $U_g\colon H_{s(g)}\to H_{r(g)}$ such that $U_{gh}=U_gU_h$  whenever $s(g)=r(h)$ and, for all Borel sections $\xi$ and $\zeta$ of $\HH$, the function $\G\ni g\mapsto (U_g\xi_{s(g)},\zeta_{r(g)})$ is Borel. One also says then that $\HH$ is a \emph{Borel $\G$-Hilbert bundle} over~$\Gu$.
\end{defn}

We remark again that in general it is more reasonable to require the functions $\G\ni g\mapsto (U_g\xi_{s(g)},\zeta_{r(g)})$ to be only locally Borel, but in this section the differences between the two requirements will be inessential.

Given a Borel action $\G\curvearrowright\HH$ and a quasi-invariant $\sigma$-finite regular Borel measure $\mu$ on $\Gu$ with Radon--Nikodym cocycle $D$, we get a representation $\pi$ of $C^*(\G)$ on $\int^\oplus_{\Gu}H_xd\mu(x)$ defined~by
\begin{equation}\label{eq:integrated1}
(\pi(f)\xi)_x:=\sum_{g\in\G^x}f(g)D(g)^{-1/2} U_g\xi_{s(g)}
\end{equation}
for $f\in C_c(\G)$ and square-integrable sections $\xi$ of $\HH$.
One calls the pair $(\mu,\G\curvearrowright\HH)$ a \emph{representation}, or a \emph{strict representation}, of $\G$ and calls $\pi$ its integrated form. As we discussed in Section~\ref{ssec:abelian}, the $C_0(\Gu)$-module $\int^\oplus_{\Gu}H_xd\mu(x)$ is topologically countably generated, hence it is topologically countably generated as  a $C^*(\G)$-module.

The following is the version of Renault's disintegration theorem~\cite{Rbook}*{Theorem~II.1.21} we are after.

\begin{thm}\label{thm:main}
Assume $\G$ is a locally compact \'etale groupoid such that the unit space $\Gu$ is Hausdorff and $\G$ can be covered by countably many open bisections. Assume $\pi\colon C^*(\G)\to B(H)$ is a  nondegenerate representation such that the $C^*(\G)$-module $H$ is topologically countably generated. 
Then there exists a quasi-invariant regular Borel probability measure $\mu$ on $\Gu$, a Borel field $\HH=(H_x)_{x\in\Gu}$ of separable Hilbert spaces over $\Gu$ and a Borel action $\G\curvearrowright\HH$ such that $H_x\ne0$ for $\mu$-a.e.~$x$ and $\pi$ is unitarily equivalent to the integrated form of $(\mu,\G\curvearrowright\HH)$.
\end{thm}

\bp
We begin by checking that the $C_0(\Gu)$-module $H$ is topologically countably generated. It suffices to show that for every $\xi\in H$ the $C_0(\Gu)$-module $\overline{\pi(C^*(\G))\xi}$ is topologically countably generated. Consider the representation $\BisG\to B(H)$, $W\mapsto u_W$, defined by Lemma~\ref{lem:bis-rep}. Choose a sequence of open bisections $W_n$ covering $\G$. Then $C_c(\G)$ is spanned by functions of the form $f*\un_{W_n}$ with $f\in C_c(r(W_n))$. Since $\pi(f*\un_{W_n})=\pi(f)u_{W_n}$ for such $f$ and $u_{W_n}\in\pi(C^*(\G))''$, it follows that the $C_0(\Gu)$-module $\overline{\pi(C^*(\G))\xi}$ is topologically generated by the vectors $u_{W_n}\xi$.

By Proposition~\ref{prop:abelian} we can then identify $H$ with $\int^\oplus_{\Gu}H_xd\mu(x)$ for a  Borel field $\HH=(H_x)_{x\in\Gu}$ of nonzero separable Hilbert spaces and a regular Borel probability measure $\mu$ in such a way that $\pi|_{C_0(\Gu)}$ becomes the representation by diagonal operators. As before, denote by $m$ the representation of $L^\infty(X,\mu)$ by diagonal operators, denote by $p_U$ the projections onto $\overline{\pi(C_0(U))H}$ for open $U$, and consider the approximate unit $\Lambda_U=\{h\in C_c(U):0\le h\le1\}$ in $C_0(U)$. Then, on the one hand, $\pi(h)\to_{h\in\Lambda_U}p_U$ strongly. On the other hand, $h\to_{h\in\Lambda_U}\un_U$ weakly$^*$ in $L^\infty(\Gu,\mu)$ by regularity of~$\mu$. By normality of the representation $m$ we conclude that $p_U=m(\un_U)$.

It follows that for every $W\in\BisG$ the partial isometry $u_W$ can be viewed as a unitary $\int^\oplus_{s(W)}H_xd\mu(x)\to\int^\oplus_{r(W)}H_xd\mu(x)$. Let $T_W\colon s(W)\to r(W)$ be the homeomorphism defined by~$W$. Since $\un_W*f=(f\circ T_W^{-1})*\un_W$ for $f\in C_c(s(W))$, the operator $u_W$ intertwines the representations $f\mapsto \pi(f)|_{p(s(W))H}$ and $f\mapsto \pi(f\circ T_W^{-1})|_{p(r(W))H}$ of $C_0(s(W))$. Since the map $(\xi_x)_{x\in r(W)}\mapsto(\xi_{T_W(x)})_{x\in s(W)}$ defines a unitary isomorphism $\int^\oplus_{r(W)}H_xd\mu(x)\cong \int^\oplus_{s(W)}H_{T_W(x)}d((T_W^{-1})_*\mu)(x)$, by Corollary~\ref{cor:abelian} we conclude that $(T_W^{-1})_*(\mu|_{r(W)})\sim\mu|_{s(W)}$. Therefore the measure $\mu$ is quasi-invariant. The same corollary shows also that $u_W$ defines an essentially unique Borel family of operators $U^W_x\colon H_x\to H_{T_W(x)}$ and $U^W_x$ is unitary for $\mu$-a.e.~$x\in s(W)$.

Consider our fixed bisections $W_n\in\BisG$ covering $\G$. Since
\begin{equation}\label{eq:consistency}
u_Vp_{s(V\cap W)}=u_{V\cap W}=u_Wp_{s(V\cap W)}
\end{equation}
for all $V,W\in\BisG$, by essential uniqueness of the families $(U^W_x)_x$ we get that $U^{V}_x=U^{W}_x$ for $\mu$-a.e.~$x\in s(V\cap W)$. Therefore by modifying inductively $(U^{W_n}_x)_x$ on sets of measure zero we may assume that $U^{W_k}_x=U^{W_n}_x$ for all~$x\in s(W_k\cap W_n)$ and $k<n$. Then we get well-defined operators $U_g\colon H_{s(g)}\to H_{r(g)}$ such that $U_g=U^{W_n}_{s(g)}$ whenever $g\in W_n$. For all Borel sections $\xi$ and $\zeta$ of $\HH$, the function $g\mapsto (U_g\xi_{s(g)},\zeta_{r(g)})$ is Borel on every bisection $W_n$, hence it is Borel on~$\G$.

Since $u_Vu_W=u_{VW}$, essential uniqueness of the families $(U^W_x)_x$ implies also that
\begin{equation}\label{eq:mult}
U^{W_k}_{T_{W_l}(x)}U^{W_l}_x=U^{W_n}_x
\end{equation}
for $\mu$-a.e.~$x\in s(W_kW_l\cap W_n)$ and all indices $k,l,n$. Let $A\subset\Gu$ be the set of all units $x$ such that either~\eqref{eq:mult} is not true for some $k,l,n$, or $U_x^{W_n}$ is not unitary for some $n$. This is a Borel set of $\mu$-measure zero. Its $\G$-saturation $[A]:=r(\G_A)=\bigcup_n T_{W_n}(A\cap s(W_n))$ is still a Borel set of measure zero by quasi-invariance of~$\mu$. We then modify $\HH$ by letting $H_x:=0$ for $x\in[A]$ and, correspondingly, put $U_g:=0$ for $g\in\G_{[A]}$ (equally well we could have taken any constant field over $[A]$ and let $U_g:=1$). Then the operators $U_g$ for $g\in\G$ define a Borel action $\G\curvearrowright\HH$.

By Corollary~\ref{cor:abelian}, the operator $u_{W}$ for $W\in\BisG$ maps $\xi\in \int^\oplus_{s(W)}H_xd\mu(x)$ into the section
$$
r(W)\ni T_W(x)\mapsto \frac{d\mu}{d((T^{-1}_W)_*\mu)}(x)^{1/2}U^W_x\xi_x=\frac{d((T_W)_*\mu)}{d\mu}(T_W(x))^{-1/2}U^W_x\xi_x.
$$
This implies that for the integrated form $\tilde\pi$ of $(\mu,\G\curvearrowright\HH)$ we have
$$
\tilde\pi(\un_{W_n}*f)=u_{W_n}\pi(f)=\pi(\un_{W_n}*f)
$$
for all $f\in C_c(s(W_n))$. Since such functions $\un_{W_n}*f$ span $C_c(\G)$, we conclude that $\tilde\pi=\pi$.
\ep

\begin{remark}
It is not difficult to see using Corollary~\ref{cor:abelian} and arguments from the above proof that the pair $(\mu,\G\curvearrowright\HH)$ is essentially uniquely determined by $\pi$, in the sense that if $(\mu',\G\curvearrowright\HH')$ is another pair with the same properties, then $\mu\sim\mu'$ and the Borel fields~$\HH$ and~$\HH'$ are $\G$-equivariantly unitarily isomorphic outside a $\G$-invariant Borel subset of $\Gu$ of $\mu$-measure zero. \ee
\end{remark}

Analyzing the proof of Theorem~\ref{thm:main} one can see that we didn't really need to define the representation of the whole inverse semigroup $\BisG$, and once this representation was constructed, to define $\G\curvearrowright\HH$ we also didn't need to remember the entire representation of $C^*(\G)$ anymore, only how the partial isometries~$u_W$ interact with functions on $\Gu$. For later use, we summarize these observations as the following corollary of the proof.

Let us first recall that an inverse subsemigroup $S\subset\BisG$ is called \emph{wide} if $\bigcup_{W\in S} W=\G$ and, for all $V,W\in S$ and $g\in V\cap W$, there exists $U\in S$ such that $g\in U\subset V\cap W$.

\begin{cor}\label{cor:inverse}
Assume $\G$ is a locally compact \'etale groupoid with Hausdorff unit space $\Gu$ and $S\subset\BisG$ is a wide inverse subsemigroup. Assume we are given a nondegenerate representation $\pi\colon C_c(\Gu)\to B(H)$ and a representation $u\colon S\to B(H)$ such that the following conditions are satisfied:
\begin{enumerate}
  \item[(i)] for every $W\in S$, the initial projection of $u_W$ is the projection $p_{s(W)}$ onto $\overline{\pi(C_c(s(W)))H}$;
  \item[(ii)] for every $W\in S$ and $f\in C_c(s(W))$, we have $u_W\pi(f)=\pi(f\circ T_W^{-1})u_W$, where we view $f\circ T_W^{-1}$ as an element of $C_c(r(W))\subset C_c(\Gu)$;
  \item[(ii)] there is a sequence $(W_n)^\infty_{n=1}$ in $S$ such that $\G=\bigcup_nW_n$;
  \item[(iv)] there is a sequence $(\xi_k)^\infty_{k=1}$ in $H$ such that the linear span of the vectors $\pi(f)u_{W}\xi_k$, with $W\in S$, $f\in C_c(r(W))$ and $k\ge1$, is dense in $H$.
\end{enumerate}
Then there exists a representation $(\mu,\G\curvearrowright\HH)$ of $\G$ such that its integrated form gives rise, up to unitary equivalence, to our given representations of~$C_c(\Gu)$ and~$S$.
\end{cor}

\bp
There are two points that have to be taken care of. First, we need to check that the conditions are enough to guarantee that $H$ is topologically countably generated as a $C_0(\Gu)$-module. For this, let us show that the linear span $H'$ of the vectors $\pi(f)u_{W_n}\xi_k$, with $f\in C_c(r(W_n))$ and $n,k\ge1$, is dense in $H$. Take $W\in S$ and $f\in C_c(r(W))$. We claim that $\pi(f)u_W\xi_k\in H'$. Since $\supp f$ is covered by the sets $r(W_n\cap W)$, the assumption that $S$ is wide and a partition of unity argument reduces the proof to the case where $\supp f\subset r(V)$ for some $n\ge1$ and $V\in S$ such that $V\subset W_n\cap W$. Then $\pi(f)=\pi(f)p_{r(V)}$ and $p_{r(V)}u_{W_n}=u_V=p_{r(V)}u_W$. Hence $\pi(f)u_W\xi_k=\pi(f)u_{W_n}\xi_k\in H'$, proving our claim.

The second point is that, given $V,W\in S$, identity~\eqref{eq:consistency} makes no sense unless $V\cap W\in S$. But similarly to the previous paragraph, we still have $u_Vp_{s(U)}=u_Wp_{s(U)}$ for any $U\in S$ such that $U\subset V\cap W$, and hence $U^V_x=U^W_x$ for $\mu$-a.e.~$x\in s(U)$. By the wideness assumption and regularity of the measure we conclude that $U^V_x=U^W_x$ for $\mu$-a.e.~$x\in s(V\cap W)$. The rest of the proof of the theorem remains unchanged.
\ep

\subsection{Representations of twisted groupoid C\texorpdfstring{$^*$}{*}-algebras}

Our next goal is to extend Theorem~\ref{thm:main} to twisted groupoid C$^*$-algebras. It will be convenient to define them in terms of Fell line bundles~\cite{MR1443836}, as then the proofs become basically the same as in the untwisted case, once we get through a fair amount of definitions.

Assume $\G$ is a locally compact \'etale groupoid with Hausdorff unit space. A \emph{Banach bundle} over $\G$ is a topological space $\LL$ together with an open continuous surjective map $p\colon\LL\to\G$ and the structure of a normed space on each fiber $\LL_g:=p^{-1}(g)$ such that
\begin{itemize}
  \item[(i)] for each $g\in\G$, the fiber $\LL_g$ is a Banach space;
  \item[(ii)] the operations of addition $\{(a,b)\in\LL\times\LL:p(a)=p(b)\}\to\LL$, $(a,b)\mapsto a+b$, and scalar multiplication $\C\times\LL\to\LL$, as well as the norm map $\LL\to[0,+\infty)$, are continuous;
  \item[(iii)] if $(a_i)_i$ is a net in $\LL$ such that $\|a_i\|\to0$ and $p(a_i)\to g$ for some $g\in\G$, then $(a_i)_i$ converges to the zero element of $\LL_g$.
\end{itemize}
Any such bundle has enough continuous sections, meaning that for every $a\in\LL$ there is a neighborhood $U$ of $p(a)$ and a continuous map $\eta\colon U\to \LL$ such that $p\circ \eta=\operatorname{id}$ and $\eta(p(a))=a$, see \cite{MR0936628}*{Appendix~C}. We will be interested only in Banach line bundles. In  this case the existence of continuous sections implies that the bundle is locally trivial. It follows that one can alternatively define a Banach line bundle as a locally trivial line bundle such that the fibers are equipped with a norm that depends continuously on the fiber, or equivalently, such that local trivializations $p^{-1}(U)\cong U\times\C$ can be chosen to be fiber-wise isometric.

A Banach line bundle $p\colon\LL\to\G$ is called a \emph{Fell line bundle} if we are in addition given antilinear maps $*\colon\LL_g\to \LL_{g^{-1}}$ and bilinear maps $\LL_g\times \LL_h\to \LL_{gh}$, $(a,b)\mapsto ab$, (when $s(g)=r(h)$) such that
\begin{itemize}
  \item[(iv)] the map $\{(a, b) \in\LL\times\LL: s(p(a)) = r(p(b))\}\to\LL$, $(a,b)\mapsto ab$, is continuous, associative and $\|ab\|=\|a\|\,\|b\|$;
  \item[(v)] the map $*\colon\LL\to\LL$ is  continuous, involutive, antimultiplicative and $\|a^*\|=\|a\|$;
  \item[(vi)] for all $a\in\LL$, we have $a^*a\ge0$.
\end{itemize}
An explanation is in order. Conditions (iv) and (v) imply that the norm satisfies $\|a^*a\|=\|a^*\|\,\|a\|=\|a\|^2$. It follows that for every unit $x\in\Gu$ the fiber $\LL_x$ is a one-dimensional C$^*$-algebra, so it can be canonically identified with $\C$. Then, for every $a\in\LL_g$, the element $a^*a$ lies in $\LL_{s(g)}=\C$ and condition (vi) requires it to be positive.

\begin{example}\label{ex:2-cocycle}
Consider the trivial Banach line bundle $\LL:=\G\times\C$ over $\G$ and a continuous normalized $\T$-valued $2$-cocycle $\omega$ on $\G$, that is, a continuous map $\omega\colon \{(g,h)\in\G\times\G: s(g)=r(h)\}\to\T$ such that $\omega(g,s(g))=\omega(r(g),g)=1$ and
$$
\omega(g,h)\omega(gh,k)=\omega(g,hk)\omega(h,k)
$$
for all composable triples $(g,h,k)$. Then we can define a Fell line bundle structure on~$\LL$~by
$$
(g,z)(h,w):=(gh,\omega(g,h)zw),\quad (g,z)^*:=(g^{-1},\overline{\omega(g^{-1},g)z}).
$$

It is not difficult to see that all Fell line bundle structures on $\LL=\G\times\C$ that are consistent with the obvious identifications of $\LL_x$ with $\C$ for $x\in\Gu$ arise this way.
\ee
\end{example}

Let $p\colon \LL\to\G$ be a Fell line bundle. Given an open subset $U\subset\G$, by a unitary section of~$\LL$ over $U$ we mean a continuous section~$\eta$ of the bundle $\LL|_U:=p^{-1}(U)$ satisfying $\|\eta(g)\|=1$ for all $g\in U$. Existence of such sections is equivalent to triviality of the line bundle $\LL|_U$. In particular, every point of $\G$ has a Hausdorff open neighborhood over which~$\LL$ admits a unitary section.

For a Hausdorff open subset $U\subset\G$, consider the space $\Gamma_c(U,\LL)$ of compactly supported continuous sections of the bundle $\LL|_U$. Similarly to how we defined $C_c(\G)$, we view the elements of $\Gamma_c(U,\LL)$ as globally defined sections by extending them by zero and denote by $\Gamma_c(\G,\LL)$ the linear span of all such spaces $\Gamma_c(U,\LL)$. This is a $*$-algebra with convolution product
\begin{equation*} \label{eprod2}
(\eta_1*\eta_2)(g) := \sum_{h \in \G^{r(g)}} \eta_1(h)\eta_2(h^{-1}g)
\end{equation*}
and involution by $\eta^*(g):=\eta(g^{-1})^*$. In view of our identifications $\LL_x=\C$ for $x\in\Gu$, we can identify $\Gamma_c(\Gu,\LL)$ with $C_c(\Gu)$. Then $C_c(\Gu)$ with its usual $*$-algebra structure becomes a $*$-subalgebra of $\Gamma_c(\G,\LL)$.

The supremum-norm on $\Gamma_c(\G,\LL)$ is defined by $\|\eta\|_\infty:=\sup_{g\in\G}\|\eta(g)\|$. The same proof as that of Lemma~\ref{lem:unibound} gives the following result.

\begin{lemma}\label{lem:unibound2}
For every $\eta\in \Gamma_c(\G,\LL)$, there is a constant $C_\eta\ge0$ depending only on $\eta$ such that if~$\pi$ is a representation of $\Gamma_c(\G,\LL)$ on a pre-Hilbert space, then $\|\pi(\eta)\|\le C_\eta$. More precisely, if $W\subset\G$ is an open bisection and $\eta\in \Gamma_c(W,\LL)$, then $\|\pi(\eta)\|\le\|\eta\|_\infty$.
\end{lemma}

We can then define a C$^*$-seminorm on $\Gamma_c(\G,\LL)$ by taking the supremum over all representations of $\Gamma_c(\G,\LL)$ by bounded operators on Hilbert spaces, or equivalently, by arbitrary operators on pre-Hilbert spaces. Similarly to the untwisted case we have regular representations $\rho^\LL_x\colon\Gamma_c(\G,\LL)\to B(\ell^2(\G_x,\LL))$ ($x\in\Gu$) defined by
\begin{equation*}\label{eq:rhox3}
(\rho^\LL_x(\eta)\xi)(g) :=
\sum_{h\in \G^{r(g)}}\eta(h) \xi(h^{-1}g),
\end{equation*}
where $\ell^2(\G_x,\LL)$ is the Hilbert space of sections $\xi$ of $\LL|_{\G_x}$ such that $\sum_{g\in\G_x}\|\xi(g)\|^2<\infty$. This implies that the C$^*$-seminorm on $\Gamma_c(\G,\LL)$ is in fact a norm dominating the supremum-norm, and if $W\subset\G$ is an open bisection and $\eta\in\Gamma_c(W,\LL)$, then
\begin{equation}\label{eq:norm-equality2}
\|\eta\|=\|\eta\|_\infty.
\end{equation}
Denote by $C^*(\G,\LL)$ the completion of $\Gamma_c(\G,\LL)$ with respect to the C$^*$-norm that we defined. When $\LL$ is the trivial Fell line bundle $\G\times\C$ (that is, the structure maps are defined as in Example~\ref{ex:2-cocycle} for $\omega\equiv1$), we get back the C$^*$-algebra $C^*(\G)$.

\begin{defn}
Assume $\HH=(H_x)_{x\in\Gu}$ is a Borel field of separable Hilbert spaces over $\Gu$.
By a \emph{Borel action} of a Fell line bundle $p\colon\LL\to\G$ on $\HH$ we mean an assignment to every $a\in\LL$ of a bounded linear operator $L(a)\colon H_{s(p(a))}\to H_{r(p(a))}$ such that it respects the linear structure, product and involution on the fibers of $\LL$ and, for all Borel sections $\xi$ and $\zeta$ of $\HH$, the function $\LL\ni a\mapsto (L(a)\xi_{s(p(a))},\zeta_{r(p(a))})$ is Borel. One also says then that $L$ is a \emph{Borel $*$-functor}.
\end{defn}

Given a Borel action $\LL\curvearrowright\HH$, for every unit $x\in\Gu$ we get a representation of $\LL_x=\C$ on~$H_x$, which must be given by a projection in $B(H_x)$. It follows that for every unitary $u\in\LL$ the operator $L(u)$ is a partial isometry. In particular, $\|L(a)\|\le\|a\|$ for all $a\in\LL$. Let us say that a Borel action $\LL\curvearrowright\HH$ is nondegenerate if the representations $L\colon\LL_x\to B(H_x)$ are unital for all $x\in\Gu$, or equivalently, if the operator $L(u)$ is unitary for every unitary $u\in\LL$.

Assume that in addition to a Borel action $\LL\curvearrowright\HH$ we are given a quasi-invariant $\sigma$-finite regular Borel measure $\mu$ on $\Gu$ with Radon--Nikodym cocycle $D$. Then we get a representation~$\pi$ of~$C^*(\G,\LL)$ on~$\int^\oplus_{\Gu}H_xd\mu(x)$ defined~by
\begin{equation}\label{eq:integrated2}
(\pi(\eta)\xi)_x:=\sum_{g\in\G^x}D(g)^{-1/2} L(\eta(g))\xi_{s(g)}
\end{equation}
for $\eta\in \Gamma_c(\G,\LL)$ and square-integrable sections $\xi$ of $\HH$. One calls the pair $(\mu,\LL\curvearrowright\HH)$ a (strict) representation of $\LL$ and calls $\pi$ its integrated form.

We are now ready to formulate and prove an extension of Theorem~\ref{thm:main} to twisted  \'etale groupoid C$^*$-algebras.

\begin{thm}\label{thm:main-twist}
Assume $\G$ is a locally compact \'etale groupoid with Hausdorff unit space~$\Gu$ and~$p\colon\LL\to\G$ is a Fell line bundle such that $\G$ can be covered by countably many open bisections over which $\LL$ admits unitary sections. Assume $\pi\colon C^*(\G,\LL)\to B(H)$ is a  nondegenerate representation such that the $C^*(\G,\LL)$-module $H$ is topologically countably generated. Then there exists a quasi-invariant regular Borel probability measure $\mu$ on $\Gu$, a Borel field $\HH=(H_x)_{x\in\Gu}$ of separable Hilbert spaces over $\Gu$ and a nondegenerate Borel action $\LL\curvearrowright\HH$ such that $H_x\ne0$ for $\mu$-a.e.~$x$ and~$\pi$ is unitarily equivalent to the integrated form of $(\mu,\LL\curvearrowright\HH)$.
\end{thm}

\bp
The proof of Theorem~\ref{thm:main} goes through with minimal changes, so we will be sketchy. Instead of the inverse semigroup $\BisG$ we consider the inverse semigroup $\BisGL$ consisting of the pairs $(W,\eta)$ such that $W$ is an open bisection of $\G$ and $\eta$ is a unitary section of $\LL$ over $W$. Note that this is indeed a semigroup, since the product of two unitary elements in~$\LL$ is unitary by multiplicativity of the norm. Similarly to Lemma~\ref{lem:bis-rep} we get a representation $\BisGL\to B(H)$, $(W,\eta)\mapsto u_{W,\eta}$, such that the initial projection of $u_{W,\eta}$ is the projection~$p_{s(W)}$ onto $\overline{\pi(C_0(s(W)))H}$ and $\pi(u_{W,\eta})\pi(f)=\pi(\eta * f)$ for $f\in C_c(s(W))$.

Fix a cover of~$\G$ by countably many open bisections $W_n$ such that there exist unitary sections $\eta_n\colon W_n\to\LL$. Since $\Gamma_c(\G,\LL)$ is spanned by the sections $f*\eta_n$ with $f\in C_c(r(W_n))$, the $C_0(\Gu)$-module $H$ is countably generated. We then get a quasi-invariant regular Borel probability measure~$\mu$ on~$\Gu$, a direct integral decomposition $H=\int^\oplus_{\Gu}H_xd\mu(x)$ and Borel families of operators $U^{W,\eta}_x\colon H_x\to H_{T_W(x)}$ that implement $u_{W,\eta}$ for $(W,\eta)\in\BisGL$ and are unitary a.e.

By essential uniqueness of the families $(U^{W,\eta}_x)_x$, we have $U^{W_k,\eta_k}_x=b_{kl}(x)U^{W_l,\eta_l}_x$ for $\mu$-a.e.~$x\in s(W_k\cap W_l)$, where $b_{kl}\colon s(W_k\cap W_l)\to\T$ is the continuous function defined by the identity $\eta_k(g)=b_{kl}(s(g))\eta_l(g)$ for $g\in W_k\cap W_l$. Similarly, we have
$
U^{W_k,\eta_k}_{T_{W_l}(x)}U^{W_l,\eta_l}_x=c^n_{kl}(x)U^{W_n,\eta_n}_x
$ 
for $\mu$-a.e.~$x\in s(W_kW_l\cap W_n)$, where $c^n_{kl}\colon s(W_kW_l\cap W_n)\to\T$ is the continuous function defined by the identity $\eta_k(g)\eta_l(h)=c^n_{kl}(s(h))\eta_n(gh)$ for $g\in W_k$ and $h\in W_l$ such that $s(g)=r(h)$ and $gh\in W_n$.

By modifying $H_x$ and the families $(U^{W_n,\eta_n}_x)_x$ on sets of measure zero we may assume that the above identities hold everywhere and the operators $U^{W_n,\eta_n}_x$ are unitary. We then get the required Borel action $\LL\curvearrowright\HH$ such that $L(a)\colon\HH_{s(p(a))}\to\HH_{r(p(a))}$ equals $c\, U^{W_n,\eta_n}_{s(p(a))}$ whenever $p(a)\in W_n$ and $a=c\,\eta_n(p(a))$ for some $c\in\C$.
\ep

\begin{remark}
Since unitary sections of  $p\colon\LL\to\G$ exist locally, the theorem applies to any $\sigma$-compact locally compact \'etale groupoid $\G$ with Hausdorff unit space. \ee
\end{remark}

\bigskip

\section{Some applications}\label{sec:applications}

\subsection{Contractivity of representations in the \texorpdfstring{$I$}{I}-norm}

Given a Fell line bundle $\LL$ over $\G$, \emph{Hahn's $I$-norm} on $\Gamma_c(\G,\LL)$ is defined by
$$
\|\eta\|_I:=\max\Big\{\sup_{x\in\Gu}\sum_{g\in\G_x}\|\eta(g)\|,\ \sup_{x\in\Gu}\sum_{g\in\G^x}\|\eta(g)\|\Big\}
$$
The following result sharpens Lemma~\ref{lem:unibound2} and allows one to equivalently define $C^*(\G,\LL)$ as the C$^*$-envelope of the involutive Banach algebra obtained by completing $\Gamma_c(\G,\LL)$ in the $I$-norm.

\begin{prop}\label{prop:I-norm}
Assume $\G$ is a locally compact \'etale groupoid with Hausdorff unit space~$\Gu$ and $p\colon\LL\to\G$ is a Fell line bundle. Then, for every representation $\pi$ of $\Gamma_c(\G,\LL)$ on a pre-Hilbert space, we have $\|\pi(\eta)\|\le\|\eta\|_I$ for all $\eta\in\Gamma_c(\G,\LL)$.
\end{prop}

We remind that by a representation on a pre-Hilbert space $H$ we mean a homomorphism into the algebra of all, possibly unbounded, linear operators on $H$ such that $(\pi(\eta)\xi,\zeta)=(\xi,\pi(\eta^*)\zeta)$ for all $\xi,\zeta\in H$.

\smallskip

Renault's disintegration theorem, in the form given in~\cite{Rbook}, implies the proposition (modulo the automatic boundedness of the operators $\pi(\eta)$) for second countable Hausdorff~$\G$ and the twists defined by continuous $2$-cocycles, see~\cite{Rbook}*{Corollary~II.1.22}. In the untwisted second countable, but not necessarily Hausdorff, case the result has been proved in~\cite{MR4750923}*{Corollary~6.6}. In the untwisted Hausdorff, but not necessarily second countable, case the result follows from the alternative approach to disintegration theory developed in~\cite{MR3851326}, see Corollary~6.2 there and Remark~\ref{rem:inductive} below. In full generality the required norm estimate has been recently obtained in~\cite{MR4948043}*{Corollary~5.8} using another newly developed disintegration framework for groupoids.  Below we reproduce the standard proof one gives once a suitable version of the disintegration theorem is available, cf.~\cite{Rbook}*{Proposition~II.1.7}.

\bp[Proof of Proposition~\ref{prop:I-norm}]
If $W_1,\dots,W_n$ are open bisections such that there exist unitary sections $\eta_i\colon W_i\to\LL$, then the groupoid $\tilde\G$ generated by $W_1,\dots,W_n$ is open in $\G$ and is covered by countably many open bisections over which $\LL$ admits unitary sections, namely, the sections obtained as convolution products of~$\eta_i$ and $\eta_i^*$. Since every section $\eta\in\Gamma_c(\G,\LL)$ lies in $\Gamma_c(\tilde\G,\LL)$ for one of such subgroupoids $\tilde\G$, it suffices to prove the proposition assuming that $\G$ itself can be covered by countably many open bisections over which $\LL$ admits unitary sections. By Lemma~\ref{lem:unibound2} it also suffices to consider (topologically) cyclic representations by bounded operators on Hilbert spaces. Given such a representation $\pi$, by Theorem~\ref{thm:main-twist} we can identify it with the integrated form of a representation $(\mu,\LL\curvearrowright\HH)$ and follow the proof of~\cite{MR0496797}*{Lemma~1.1}. Namely, take $\xi,\zeta\in H=\int^\oplus_{\Gu}H_xd\mu(x)$. Using the Cauchy--Schwarz inequality thrice we get for all $\eta\in\Gamma_c(\G,\LL)$ that
\begin{align*}
|(\pi(\eta)\xi,\zeta)| &=\Big|\int_{\Gu}\sum_{g\in\G^x}\big(L(\eta(g))\xi_{s(g)},\zeta_x\big)D(g)^{-1/2}d\mu(x)\Big| \\
   & \le \int_{\Gu}\sum_{g\in\G^x}\|\eta(g)\|\,\|\xi_{s(g)}\|\,\|\zeta_x\|\,D(g)^{-1/2}d\mu(x)\\
   & \le \int_{\Gu}\Big(\sum_{g\in\G^x}\|\eta(g)\|\,\|\xi_{s(g)}\|^2D(g)^{-1}\Big)^{1/2}\Big(\sum_{h\in\G^x}\|\eta(h)\|\,\|\zeta_x\|^2\Big)^{1/2}d\mu(x) \\
   & \le \Big(\int_{\Gu}\sum_{g\in\G^x}\|\eta(g)\|\,\|\xi_{s(g)}\|^2D(g)^{-1}d\mu(x)\Big)^{1/2}\Big(\int_{\Gu}\sum_{h\in\G^x}\|\eta(h)\|\,\|\zeta_x\|^2d\mu(x)\Big)^{1/2} \\
   &= \Big(\int_{\Gu}\sum_{g\in\G_x}\|\eta(g)\|\,\|\xi_{x}\|^2d\mu(x)\Big)^{1/2}\Big(\int_{\Gu}\sum_{h\in\G^x}\|\eta(h)\|\,\|\zeta_x\|^2d\mu(x)\Big)^{1/2}\\
   &\le\|\eta\|_I\|\xi\|\,\|\zeta\|,
\end{align*}
where we used~\eqref{eq:RN} in the next to last line. This proves the proposition.
\ep

\begin{remark}\label{rem:inductive}
When $\G$ is Hausdorff, the inductive limit topology on $C_c(\G)$ is defined as the strongest topology making the embeddings $C_c(U)\hookrightarrow C_c(\G)$ continuous for all relatively compact open sets $U\subset \G$, where the spaces $C_c(U)$ are equipped with the topology defined by the supremum-norm. Since the norms $\|\cdot\|_\infty$ and $\|\cdot\|_I$ are equivalent on every such space~$C_c(U)$, the inductive limit topology is stronger than the $I$-norm topology. It follows that every representation $\pi\colon C_c(\G)\to B(H)$ is continuous in the inductive limit topology. This, however, does not require much effort and can be quickly deduced already from~\eqref{eq:upper-bound}, see the first part of the proof of~\cite{SSW}*{Lemma~3.2.4}.

The inductive limit and the $I$-norm topologies are in general different.\footnote{Contrary to the claim in the second part of the proof of~\cite{SSW}*{Lemma~3.2.4} that the former is weaker than the latter; see also \url{https://www.aidansims.com/files/GroupoidsCRMerrata.pdf}.} In order to see this, take any Hausdorff $\sigma$-compact, but noncompact, locally compact \'etale groupoid $\G$ and choose an increasing sequence of open subsets $U_n\subset\G$ such that $\bar U_n$ is compact, $\bar U_n\subsetneq U_{n+1}$ and $\bigcup_nU_n=\G$. Using inductively the Hahn--Banach theorem one can easily construct a linear functional~$\varphi$ on~$C_c(\G)$ such that its restriction to $(C_c(U_n),\|\cdot\|_I)$ is bounded, but has norm~$\ge n$. Then~$\varphi$ is continuous in the inductive limit topology, but is unbounded with respect to the $I$-norm. In fact, by elaborating on this argument it is not difficult to see that the inductive limit topology on $C_c(\G)$ is not even first countable, so it is not metrizable. \ee
\end{remark}

Let us also record the following automatic boundedness result, which relies only on Lemma~\ref{lem:unibound2}.

\begin{prop}\label{prop:pos}
Assume $\G$ is a locally compact \'etale groupoid with Hausdorff unit space~$\Gu$ and $p\colon\LL\to\G$ is a Fell line bundle. Then a positive linear functional $\varphi$ on $\Gamma_c(\G,\LL)$ is bounded with respect to the norm on $C^*(\G,\LL)$ if and only if its restriction to $C_c(\Gu)$ is bounded.
\end{prop}

\bp
The ``only if'' direction is clear. For the ``if'' direction, assume $\varphi$ is a positive linear functional on $\A:=\Gamma_c(\G,\LL)$ such that $C:=\|\varphi|_{C_c(\Gu)}\|$ is finite. Consider the GNS-map $\Lambda\colon\A\to H$ associated to $\varphi$, so $H$ is a Hilbert space and $\Lambda$ is a linear map with dense image such that $(\Lambda(a),\Lambda(b))=\varphi(b^*a)$. The action of $\A$ on itself by multiplication defines a representation of~$\A$ on the pre-Hilbert space $\Lambda(\A)$. This representation is contractive by the definition of the norm on $C^*(\G,\LL)$, hence it extends to a representation $\pi\colon C^*(\G,\LL)\to B(H)$. Take $h\in C_c(\Gu)$, $0\le h\le1$. Then $(\pi(\cdot)\Lambda(h),\Lambda(h))$ is a positive linear functional of norm not greater than $\|\Lambda(h)\|^2=\varphi(h^2)\le C$. Thus, for every $\eta\in \Gamma_c(\G,\LL)$ we have $|\varphi(h*\eta*h)|\le C\|\eta\|$. Since for any given $\eta$ we can find $h\in C_c(\Gu)$ such that $0\le h\le 1$ and $\eta=h*\eta=\eta*h$, we get the result.
\ep

We finish this subsection with a brief discussion of \emph{ample} groupoids $\G$, which by definition means that the unit space $\Gu$ is totally disconnected. Denote by $A(\G)\subset C_c(\G)$ the linear span of characteristic functions $\un_W$ of compact open bisections $W$. It is a $*$-subalgebra of the convolution algebra $C_c(\G)$, called the (complex) \emph{Steinberg algebra} of~$\G$.

The following proposition generalizes a result of Clark and Zimmerman~\cite{MR4750923}*{Theorem~7.1}, proved for second countable groupoids.

\begin{prop}\label{prop:Steinberg}
Assume $\G$ is an ample locally compact \'etale groupoid with Hausdorff unit space~$\Gu$. Then, for every representation $\pi$ of the Steinberg algebra $A(\G)$ on a pre-Hilbert space, we have $\|\pi(a)\|\le\|a\|_I$ for all $a\in A(\G)$. In particular, the universal C$^*$-completion of~$A(\G)$ is well-defined, and then the embedding $A(\G)\hookrightarrow C_c(\G)$ extends to an isomorphism of this completion onto $C^*(\G)$.
\end{prop}

\bp
Since the operator $\pi(\un_W)$ is a partial isometry for every compact open bisection $W$, we have $\|\pi(\un_W)\|\le1$. Hence $A(\G)$ has a universal C$^*$-completion $\overline{A(\G)}$, and any representation of~$A(\G)$ on a pre-Hilbert space $H$ extends to a representation of $\overline{A(\G)}$ onto the completion~$\bar H$ of~$H$. Since the subalgebra $A(\Gu)\subset A(\G)$ of compactly supported locally constant functions on~$\Gu$ is an increasing union of finite dimensional $*$-algebras, the supremum-norm on it is the unique C$^*$-norm, so the closure of $A(\Gu)$ in $\overline{A(\G)}$ can be identified with $C_0(\Gu)$.

It is clear that $A(\G)$ is dense in $C^*(\G)$, since $A(\Gu)$ is dense in $C_c(\Gu)$ and $C_c(\G)$ is spanned by functions of the form $f*\un_W$, where $W$ is a compact open bisection and $f\in C_c(r(W))=C(r(W))$. In view of Proposition~\ref{prop:I-norm} all that remains to show to finish the proof is that every Hilbert space representation $\pi\colon A(\G)\to B(H)$ extends (necessarily uniquely) to a representation of $C_c(\G)$. Similarly to the proof of Proposition~\ref{prop:I-norm}, it is enough to be able to extend~$\pi|_{A(\tilde \G)}$ to~$C_c(\tilde \G)$ for every open subgroupoid $\tilde \G\subset\G$ generated by finitely many compact open bisections. It follows that without loss of generality we may assume that $\G$ can be covered by countably many compact open bisections $W_n$, and then we may assume in addition that $\pi$ is (topologically) cyclic.

The representation $\pi|_{A(\Gu)}$ extends by continuity to a representation of $C_0(\Gu)$. We continue to write~$\pi(f)$ ($f\in C_0(\Gu)$) for the operators of this representation. Denote by $\BiscG\subset\BisG$ the inverse subsemigroup of compact open bisections. Clearly, it is wide. We get a representation of~$\BiscG$ on~$H$ defined by $u_W:=\pi(\un_W)$. We have $u_W\pi(f)=\pi(f\circ T_W^{-1})u_W$ for all $f\in A(s(W))$, hence for all $f\in C(s(W))$, where $T_W\colon s(W)\to r(W)$ is the homeomorphism defined by $W$ and $f\circ T_W^{-1}$ is viewed as an element of $C(r(W))\subset C_0(\Gu)$.

At this point we can apply Corollary~\ref{cor:inverse}, with $S:=\BiscG$, and conclude that the representation~$\pi$ of~$A(\G)$ is the integrated form of a representation $(\mu,\G\curvearrowright\HH)$ of $\G$. The latter representation can then be integrated to a representation of $C_c(\G)$, giving the required extension of~$\pi$ (which in the end is simply given by mapping $f*\un_W\in C_c(\G)$ into $\pi(f)u_W$ for $W\in\BiscG$ and $f\in C(r(W))$).
\ep

It is, in fact, possible to prove the isomorphism $\overline{A(\G)}\cong C^*(\G)$ without using any disintegration results, see Remark~\ref{rem:BKM2}.

\subsection{Crossed products by inverse semigroups}\label{ssec:inverse}

Assume $X$ is a Hausdorff locally compact space and $S$ is an inverse semigroup. By an action of $S$ on $X$ one means a homomorphism of $S$ into the inverse semigroup of partial homeomorphisms of $X$ (that is, homeomorphisms between open subsets of $X$) such that the domains of definition of the partial homeomorphisms corresponding to $s\in S$ cover $X$. Denote by $E(S)$ the set of idempotents in $S$. For every $e\in E(S)$, denote by $D_e\subset X$ the corresponding domain of definition, so that $e$ acts by the identity map on~$D_e$. Then $s\in S$ acts by a homeomorphism $D_{s^{-1}s}\to D_{ss^{-1}}$.

To any action $S\curvearrowright X$ one can associate a C$^*$-algebra $C_0(X)\rtimes S$ and an \'etale groupoid $S\ltimes X$. Let us start with the latter, see~\cite{MR1724106}*{Section~3.3}, \cite{MR2419901}*{Section~4} for details. Denote by $\Omega\subset S\times X$ the set of pairs $(s,x)$ such that $x\in D_{s^{-1}s}$. Define an equivalence relation on $\Omega$ by
$$
(s,x)\sim (t,y)\quad\text{iff}\quad x=y\ \ \text{and}\ \ se=te\ \ \text{for some}\ \ e\in E(S)\ \ \text{with}\ \ x\in D_e.
$$
Denote by $[s,x]$ the equivalence class of $(s,x)\in\Omega$. The quotient set $S\ltimes X:=\Omega/_\sim$ is a groupoid with product
$$
[s,tx][t,x]:=[st,x].
$$
Its unit space is identified with $X$ by identifying $x\in X$ with $[e,x]$ for any $e\in E(S)$ satisfying $x\in D_e$. Therefore the source and range maps are $[s,x]\mapsto x$ and $[s,x]\mapsto sx$, resp., and the inverse is given by $[s,x]^{-1}=[s^{-1},sx]$.

As a basis of topology on $S\ltimes X$ one takes the sets $\{[s,x]:x\in U\}$ for $s\in S$ and open $U\subset D_{s^{-1}s}$. Then $S\ltimes X$ becomes a locally compact, but not necessarily Hausdorff, \'etale groupoid. Note that its unit space $X$ carries the original topology. Observe also that every $s\in S$ defines an open bisection
\begin{equation}\label{eq:Ws}
W_s:=\{[s,x]:x\in D_{s^{-1}s}\}\subset S\ltimes X.
\end{equation}
The map $S\to\Bis(S\ltimes X)$, $s\mapsto W_s$, is a homomorphism of inverse semigroups. In general it is neither injective nor surjective, but its image is always a wide inverse subsemigroup of $\Bis(S\ltimes X)$.

\smallskip

Let us turn to the construction of $C_0(X)\rtimes S$, see~\cite{MR1456588}, \cite{MR2419901}*{Section~9} for details and \cite{MR2881538} for a Fell bundle generalization. Consider the linear space $\C_X[S]$ of finite formal combinations $\sum_{s\in S} f_s s$ of elements of $S$ with coefficients $f_s\in C_c(D_{ss^{-1}})$. This is a $*$-algebra with product
$$
(f s)(g\, t):= (fg(s^{-1}\cdot))st
$$
and involution $(fs)^*:=\overline{f(s\,\cdot)}s^{-1}$. A pair of representations $\pi\colon C_c(X)\to B(H)$ and $v\colon S\to B(H)$ is called \emph{covariant} if
\begin{itemize}
  \item[(i)] for every $e\in E(S)$, the operator $v_e$ is the projection onto $\overline{\pi(C_c(D_e))H}$;
  \item[(ii)] for all $s\in S$ and $f\in C_c(D_{s^{-1}s})$, we have $v_s\pi(f)=\pi(f(s^{-1}\cdot))v_s$, where $f(s^{-1}\cdot)$ is viewed as an element of $C_c(D_{ss^{-1}})$.
\end{itemize}
Every such pair defines a representation $\rho\colon \C_X[S]\to B(H)$ by $\rho(f s):=\pi(f)v_s$ for $s\in S$ and $f\in C_c(D_{ss^{-1}})$. In general not all representations of $\C_X[S]$ arise this way, only so called \emph{admissible} ones. Taking the supremum of $\|\rho(a)\|$ for $a\in\C_X[S]$ over all admissible~$\rho$, we get a C$^*$-seminorm on $\C_X[S]$. By factoring out its kernel and completing we obtain a C$^*$-algebra $C_0(X)\rtimes S$.

\smallskip

We now have two C$^*$-algebras associated with an action $S\curvearrowright X$, the crossed product $C_0(X)\rtimes S$ and the groupoid C$^*$-algebra $C^*(S\ltimes X)$. To compare them, observe that there is a canonical $*$-homomorphism $\theta\colon\C_X[S]\to C_c(S\ltimes X)$ defined by
$$
\sum_{s\in S} f_s s\mapsto \sum_{s\in S} f_s*\un_{W_s},
$$
where $W_s$ is given by~\eqref{eq:Ws}. Since the bisections $W_s$ cover $S\ltimes X$, it is clear that $\theta$ is surjective. By Lemma~\ref{lem:bis-rep} any representation $\rho \colon C^*(S\ltimes X)\to B(H)$ gives rise to a representation of $\Bis(S\ltimes X)$, hence of $S$. This shows that $\rho\circ\theta$ is an admissible representation of $\C_X[S]$. It follows that $\theta$ extends to a surjective homomorphism $C_0(X)\rtimes S\to C^*(S\ltimes X)$.

The following theorem generalizes a result of Paterson~\cite{MR1724106}*{Corollary~3.3.2}.

\begin{thm}\label{thm:inverse1}
For any action $S\curvearrowright X$ of an inverse semigroup $S$ on a Hausdorff locally compact space $X$, the canonical map $C_0(X)\rtimes S\to C^*(S\ltimes X)$ is an isomorphism.
\end{thm}

The result has been known for a while at least when $S\ltimes X$ is either second countable or Hausdorff~\cite{MR2419901}*{Theorem 9.8}, \cite{MR2881538}*{Theorem~2.15}, \cite{MR3851326}*{Theorem~7.6}. It is stated (even in a more general form) in \cite{MR3619758}*{Theorem~5.4}, but it seems the proof given there requires either a suitable disintegration theorem or an alternative definition of the algebras involved. It can in any case be quickly deduced from~\cite{MR4948043}*{Corollaries~4.26 and~4.27}, see also Remark~\ref{rem:BKM} below. We will give a short proof based on Corollary~\ref{cor:inverse}.

\bp[Proof of Theorem~\ref{thm:inverse1}]
Since the canonical map $C_0(X)\rtimes S\to C^*(S\ltimes X)$ is surjective, in order to prove the theorem it is enough to show that every representation $\rho\colon C_0(X)\rtimes S\to B(H)$ defined by a covariant pair $(\pi,v)$ factors through $C^*(S\ltimes X)$. Since the map $\theta\colon\C_X[S]\to C_c(S\ltimes X)$ is surjective, for this, in turn, it suffices to show that the restriction of $\rho$ to $\C_X[S]$ factors through~$C_c(S\ltimes X)$.

Take $s_1,\dots,s_n\in S$ and $\xi_1,\dots,\xi_n\in H$. Denote by $\tilde\G$ the open subgroupoid of $S\ltimes X$ generated by the bisections $W_{s_k}$, $1\le k\le n$, and let $\tilde X\subset X$ be its unit space. Let $\tilde S\subset S$ be the inverse subsemigroup consisting of elements $s\in S$ such that $W_s\subset\tilde\G$. Then $\tilde \G$ can be identified with $\tilde S\ltimes\tilde X$. Consider the closure $\tilde H\subset H$ of the subspace spanned by $\rho(\C_{\tilde X}[\tilde S])\xi_k$, $1\le k\le n$, and the representation $\tilde\rho(a):=\rho(a)|_{\tilde H}$ of $\C_{\tilde X}[\tilde S]$ on $\tilde H$. Then it suffices to show that $\tilde\rho$ factors through $C_c(\tilde S\ltimes\tilde X)$, for all possible choices of $s_k$ and $\xi_k$. It follows that without loss of generality we may assume that $S\ltimes X$ is covered by countably many bisections $W_s$ and there is a finite number of vectors $\xi_k\in H$ such that the linear span of $\rho(\C_X[S])\xi_k$ is dense in $H$.

We claim that the representation $v\colon S\to B(H)$ factors through the image of~$S$ in~$\Bis(S\ltimes X)$ under the map $s\mapsto W_s$. Indeed, assume $W_s=W_t$ for some $s,t\in S$. Then $D_{s^{-1}s}=D_{t^{-1}t}$. In order to show that $v_s=v_t$, by condition (i) in the definition of covariant pairs it suffices to check that $v_s\pi(f)=v_t\pi(f)$ for all $f\in C_c(D_{s^{-1}s})$. Given such a function $f$, for every point $x\in\supp f$, we can find an idempotent $e_x\in S$ such that $x\in D_{e_x}$ and $se_x=te_x$. Since $\supp f$ is compact, a partition of unity argument then reduces the problem to showing that $v_s\pi(f)=v_t\pi(f)$ under the assumption that there exists an idempotent $e\in S$ satisfying $\supp f\subset D_e\cap D_{s^{-1}s}$ and $se=te$. But then the equality is obvious, since $\pi(f)=v_e\pi(f)$ and $v_s v_e=v_{se}=v_{te}=v_tv_e$. Thus, we get a representation $u$ of the wide inverse subsemigroup $\{W_s: s\in S\}\subset \Bis(S\ltimes X)$ defined unambiguously by $u_{W_s}:=v_s$.

We are now in a position to apply Corollary~\ref{cor:inverse} and conclude that the pair~$(\pi,u)$ arises from a representation $(\mu,(S\ltimes X)\curvearrowright\HH)$. By integrating we then get a representation $C_c(S\ltimes X)\to B(H)$, explicitly given~by
$$
f* \un_{W_s}\mapsto\pi(f)u_{W_s}=\pi(f)v_s=\rho(f s)
$$
for $s\in S$ and $f\in C_c(D_{ss^{-1}})$. This shows that $\rho$ indeed factors through $C_c(S\ltimes X)$ and completes the proof of the theorem.
\ep

\begin{remark}\label{rem:BKM}
Once we checked that $v\colon S\to B(H)$ factors through the image of $S$ in $\Bis(S\ltimes X)$, the desired conclusion that $\rho\colon\C_X[S]\to B(H)$ factors through $C_c(S\ltimes X)$ was an immediate application of the disintegration theorem in the form of Corollary~\ref{cor:inverse}. Although this argument is effortless, it is in fact possible to prove this factorization by a clever direct argument, no disintegration is needed, see the proof of~\cite{MR4948043}*{Proposition~4.23}. \ee
\end{remark}

Our next goal is to obtain a crossed product decomposition of $C^*(\G)$ for an \'etale groupoid $\G$. Let us first make a few remarks about homomorphisms between groupoid C$^*$-algebras.

Assume $\G_1$ and $\G_2$ are locally compact \'etale groupoids with Hausdorff unit spaces, and $\pi\colon\G_1\to\G_2$ is a continuous groupoid homomorphism that restricts to a homeomorphism of~$\G^{(0)}_1$ onto~$\G^{(0)}_2$. Then $\pi$ is a local homeomorphism. Moreover, for any open bisection $W\subset\G_1$, the image $\pi(W)$ is an open bisection of $\G_2$ and $\pi$ defines a homeomorphism $W\to\pi(W)$. It follows that $\pi$ defines a homomorphism $\Bis(\G_1)\to\Bis(\G_2)$, and then we get a $*$-homomorphism $\pi_*\colon C_c(\G_1)\to C_c(\G_2)$ that maps $f\in C_c(W)$ for an open bisection $W\subset\G_1$ into $f\circ\pi^{-1}_W\in C_c(\pi(W))$, where $\pi_W^{-1}\colon \pi(W)\to W$ is the inverse of $\pi|_W$. In other words, $\pi_*$ is given by
$$
\pi_*(f)(g)=\sum_{h\in\pi^{-1}(g)}f(h)\quad\text{for}\quad f\in C_c(\G_1).
$$
By passing to completions we get a homomorphism $C^*(\G_1)\to C^*(\G_2)$, which we continue to denote by~$\pi_*$.

\begin{lemma}\label{lem:iso}
The map $\pi_*\colon C^*(\G_1)\to C^*(\G_2)$ is an isomorphism if and only if $\pi$ is bijective, in which case $\pi$ is a topological isomorphism.
\end{lemma}

\bp
The ``if'' direction is obvious. For the opposite direction, assume $\pi_*$ is an isomorphism. Assume first that $\pi$ is not surjective and take $g\in\G_2\setminus\pi(\G_1)$. Consider the representation~$\rho_{s(g)}$ defined by~\eqref{eq:rhox2} and the bounded linear functional $\varphi:=(\rho_{s(g)}(\cdot)\delta_{s(g)},\delta_g)$ on $C^*(\G_2)$. Then $\varphi(f)=f(g)$ for $f\in C_c(\G_2)$. This implies that $\varphi\ne0$ and at the same time $\varphi$ vanishes on~$\pi_*(C_c(\G_1))$, hence on $\pi_*(C^*(\G_1))$. This contradiction proves that $\pi$ must be surjective.

Assume now that $\pi$ is not injective. Then there exist $x\in\G^{(0)}_1$ and $g\in(\G_1)^x_x\setminus\{x\}$ such that $\pi(g)=\pi(x)$. Put $y:=\pi(x)$. Take any open bisection $W\subset\G_1$ containing $g$. Let $U:=s(W)$, choose any function $h\in C_c(U)$ such that $h(x)=1$ and consider $f:=\un_W*h-h$. Then $f|_{(\G_1)^x_x}=\delta_g-\delta_x\ne0$, while $\pi_*(f)|_{(\G_2)^y_y}=0$. To see that this leads to a contradiction, consider a net of functions $q_i\in C_c(U)$ such that $q_i(x)=1$ for all $i$ and for any neighborhood $V$ of $x$ we eventually have $\supp q_i\subset V$. Then, by \cite{MR4742724}*{Theorem~1.8}, we have
$$
\lim_i\|q_i*f*q_i\|_{C^*(\G_1)}=\|f|_{(\G_1)^x_x}\|_{C^*((\G_1)^x_x)}\ne0
$$
and $\lim_i\|\pi_*(q_i*f*q_i)\|_{C^*(\G_2)}=\|\pi_*(f)|_{(\G_2)^y_y}\|_{C^*((\G_2)^y_y)}=0$, which contradicts the fact that~$\pi_*$ is isometric by assumption. Hence $\pi$ is injective.
Then, since $\pi$ is a bijective local homeomorphism, it must be a homeomorphism.
\ep

Assume now that we are given  a locally compact \'etale groupoid $\G$ with Hausdorff unit space and an inverse subsemigroup $S\subset\BisG$. Assuming that $\bigcup_{W\in S}r(W)=\Gu$, we get an action $S\curvearrowright\Gu$ such that each $W\in S$ acts by the usual homeomorphism $s(W)\to r(W)$ defined by~$W$. Consider the corresponding groupoid $S\ltimes\Gu$. We have a canonical continuous groupoid homomorphism $\pi_S\colon S\ltimes\Gu\to\G$ that maps $[s,x]$ into the unique element of $s\cap\G_x$. Hence we get a homomorphism $(\pi_S)_*\colon C^*(S\ltimes\Gu)\to C^*(\G)$. Composing it with the canonical homomorphism $C_0(\Gu)\rtimes S\to C^*(S\ltimes\Gu)$, we get a homomorphism $C_0(\Gu)\rtimes S\to C^*(\G)$, explicitly given~by
$$
\sum_{s\in S}f_ss\mapsto\sum_{s\in S} f_s*\un_s.
$$

The following theorem generalizes and refines another (but closely related) result of Paterson~\cite{MR1724106}*{Theorem~3.3.1}, see also \cite{MR1671944}*{Theorem~8.1}, \cite{MR2419901}*{Proposition~9.9}, \cite{MR3619758}*{Corollary~5.6}, \cite{MR3851326}*{Theorem~7.6}, \cite{MR4948043}*{Corollary~4.27}.

\begin{thm}\label{thm:inverse2}
Assume $\G$ is a locally compact \'etale groupoid with Hausdorff unit space~$\Gu$ and $S\subset\BisG$ is an inverse subsemigroup such that $\Gu=\bigcup_{W\in S}r(W)$. Then the canonical map $C_0(\Gu)\rtimes S\to C^*(\G)$ is an isomorphism if and only if the canonical map $\pi_S\colon S\ltimes\Gu\to\G$ is a bijection (hence a topological groupoid isomorphism).
\end{thm}

\bp
As $C_0(\Gu)\rtimes S\cong C^*(S\ltimes\Gu)$ by Theorem~\ref{thm:inverse1}, what we equivalently need to prove is that $\pi_S\colon S\ltimes\Gu\to\G$ induces an isomorphism of groupoid C$^*$-algebras if and only if it is bijective. This is true by Lemma~\ref{lem:iso}.
\ep

It is not difficult to check that bijectivity of $\pi_S\colon S\ltimes\Gu\to\G$ is equivalent to wideness of $S\subset\BisG$, see \cite{MR2419901}*{Proposition~5.4}, \cite{MR4246403}*{Proposition~2.2}.

\begin{remark}\label{rem:BKM2}
As was indicated in Remark~\ref{rem:BKM}, it is possible to prove Theorem~\ref{thm:inverse1}, and hence Theorem~\ref{thm:inverse2}, without using any disintegration results. One can then use the decomposition $C^*(\G)= C_0(\Gu)\rtimes S$ for some wide inverse subsemigroup $S\subset\BisG$ to systematically study the representation theory of $C^*(\G)$. This is the approach taken in~\cite{MR4948043}. For example, looking at the proof of Proposition~\ref{prop:Steinberg} from this angle, we see that any representation~$\pi$ of~$A(\G)$ defines a covariant pair of representations for the action $\BiscG\curvearrowright\Gu$, and conclude that this is already enough to extend $\pi$ to~$C_c(\G)$.\ee
\end{remark}

\subsection{KMS-states and weights}\label{ssec:KMS}

Consider again a Fell line bundle $\LL$ over an \'etale groupoid~$\G$. Assume $c\colon\G\to\R$ is a continuous homomorphism, also called a $1$-cocycle. Then we get a one-parameter group of automorphisms $\sigma^c_t$ of $C^*(\G,\LL)$ defined by
\begin{equation}\label{eq:dynamics}
\sigma^c_t(\eta)(g):=e^{it c(g)}\eta(g)
\end{equation}
for $\eta\in\Gamma_c(\G,\LL)$ and $g\in\G$. In this subsection we discuss the structure of $\sigma^c$-KMS-states and weights on $C^*(\G,\LL)$.

\smallskip

We refer the reader to \cites{Ku,MR4864224} for a thorough introduction to KMS-weights, let us only briefly recall some of the key definitions and facts. A weight on a C$^*$-algebra $A$ is a map $\psi\colon A_+\to[0,+\infty]$ that is additive and positively homogeneous. For any weight $\psi$, define
$$
\mathcal{N}_\psi:=\{a\in A:\psi(a^*a)<\infty\},\quad\mathcal{M}_\psi:=\mathcal{N}_\psi^*\mathcal{N}_\psi=\operatorname{span}\{a^*b: a,b\in \mathcal{N}_\psi\}.
$$
Then $\mathcal{M}_\psi$ is a $*$-subalgebra of $A$ that coincides with $\operatorname{span}\{a\in A_+:\psi(a)<\infty\}$, and $\psi$ extends uniquely to a positive linear functional on~it.

Assume $\sigma=(\sigma_t)_{t\in\R}$ is a one-parameter group of automorphisms of~$A$. An element $a\in A$ is called \emph{$\sigma$-analytic} if the map $\R\to A$, $t\mapsto\sigma_t(a)$, has an analytic extension to $\C$. Given $\beta\in\R$, a weight $\psi$ on $A$ is called a \emph{$\sigma$-KMS$_\beta$-weight} if it is $\sigma$-invariant, lower semicontinuous, densely defined and
\begin{equation}\label{eq:KMS}
\psi(a^*a)=\psi\big(\sigma_{-\frac{i\beta}{2}}(a)\sigma_{-\frac{i\beta}{2}}(a)^*\big)
\end{equation}
for all $\sigma$-analytic elements $a\in A$.

When $\sigma$ is trivial, so that we deal with lower semicontinuous densely defined traces on $A$, by a classical result of Pedersen~\cite{MR0212582} there is a one-to-one correspondence between such traces and the positive tracial linear functionals on the Pedersen ideal $\Ped(A)$ of $A$. The following proposition gives a partial extension of this correspondence to KMS-weights.

\begin{prop}\label{prop:Ped}
Assume $A$ is a C$^*$-algebra and $\sigma=(\sigma_t)_{t\in\R}$ is a one-parameter group of automorphisms of~$A$.  Assume $\A\subset A$ is a dense $*$-subalgebra satisfying the following properties:
\begin{itemize}
  \item[(a)] $\A$ consists of $\sigma$-analytic elements and $\sigma_z(\A)\subset\A$ for all $z\in\C$;
  \item[(b)] for $\BB:=\A\cap\Ped(A^\sigma)$, we have $\A\subset \BB A\BB$;
  \item[(c)] $\A=\A^2$;
  \item[(d)] every representation of~$\A$ on a pre-Hilbert space is contractive.
\end{itemize}
Then, for every $\beta\in\R$, the restriction map $\psi\mapsto\psi|_\A$ defines a one-to-one correspondence between the $\sigma$-KMS$_\beta$-weights $\psi$ on $A$ and the positive linear functionals $\varphi$ on $\A$ satisfying the following properties:
\begin{enumerate}
  \item[(i)] $\varphi(\sigma_z(a))=\varphi(a)$ for all $a\in\A$ and $z\in\C$;
  \item[(ii)] the function $\C\ni z\mapsto\varphi(a\sigma_z(b))$ is analytic for all $a,b\in\A$;
  \item[(iii)] $\varphi(ab)=\varphi(b\sigma_{i\beta}(a))$ for all $a,b\in\A$.
\end{enumerate}
\end{prop}

Here $A^\sigma$ denotes the fixed-point subalgebra of $A$ and the positivity means as usual that $\varphi(a^*a)\ge0$ for all $a\in\A$. 
When $\sigma$ is trivial, we may take~$\A$ to be any dense $*$-subalgebra of~$\Ped(A)$ such that its positive cone $\A_+:=\A\cap A_+$ is closed under taking square roots and spans~$\A$; for example, the Pedersen ideal $\Ped(A)$ itself has these properties, see~\cite{MR3839621}*{Section~5.6}. Indeed, then conditions (a)--(c) are obviously satisfied, while (d) is verified using the same arguments as in the proof of Lemma~\ref{lem:unibound}. For general $\sigma$, to have such an algebra $\A$ at the very least we need~$A^\sigma$ to be a nondegenerate C$^*$-subalgebra of $A$, but it is not clear to us whether this condition is enough.

\bp[Proof of Proposition~\ref{prop:Ped}]
If $\psi$ is a $\sigma$-KMS$_\beta$-weight, then $\Ped(A^\sigma)\subset\M_\psi$ by \cite{MR4592883}*{Proposition~3.1}. Since $\NN_\psi\subset A$ is a left ideal, we have $\M_\psi A\M_\psi\subset\NN_\psi^* A\NN_\psi\subset\M_\psi$. Condition (b) implies then that $\A\subset\M_\psi$, so the restriction $\varphi:=\psi|_\A$ is a well-defined positive linear functional. The same condition and density of $\A$ in $A$ imply that~$\bar\BB$ is a nondegenerate C$^*$-subalgebra of $A$. Let $(e_i)_i$ be an approximate unit in $A$ consisting of elements of $\BB$. By \cite{MR2056837}*{Lemma~3.1(iii),(iv)}, we have $\psi(e_iae_i)\to\psi(a)$ for all $a\in A_+$. Since the positive linear functionals $\psi(e_i\cdot e_i)$ on $A$ are determined by their restrictions to $\A$, we thus see that $\psi$ is determined by $\varphi$. Properties (i)--(iii) are standard properties of KMS-weights valid for any $\sigma_z$-invariant subalgebra of $\M_\psi$.

To finish the proof we need to show that any positive linear functional $\varphi$ on $\A$ satisfying conditions (i)--(iii) can be extended to a $\sigma$-KMS$_\beta$-weight on $A$. As in the proof of Proposition~\ref{prop:pos}, consider the GNS-map $\Lambda\colon\A\to H$ associated to $\varphi$. By condition (d) we get a representation $\pi\colon A\to B(H)$ uniquely determined by $\pi(a)\Lambda(b)=\Lambda(ab)$ for $a,b\in\A$.

Since the kernel of $\Lambda\colon\A\to H$ consists of the elements $a\in\A$ such that $\varphi(ba)=0$ for all $b\in\A$, by property (i) this kernel is $\sigma_z$-invariant for all $z\in\C$. Since
$$
\varphi(a^*a)=\varphi\big(\sigma_{-\frac{i\beta}{2}}(a)\sigma_{\frac{i\beta}{2}}(a^*)\big)=\varphi\big(\sigma_{\frac{i\beta}{2}}(a^*)^*\sigma_{\frac{i\beta}{2}}(a^*)\big)
$$
by~(i) and~(iii) and the property $\sigma_z(a)^*=\sigma_{\bar z}(a^*)$, we see that if $a\in\ker\Lambda$, then $\sigma_{\frac{i\beta}{2}}(a^*)\in\ker\Lambda$, hence $a^*\in\ker\Lambda$. Therefore $\ker\Lambda$ is a two-sided $*$-ideal in $\A$. It follows that $\UU:=\Lambda(\A)$ can be given the structure of a $*$-algebra with product $\Lambda(a)\Lambda(b):=\Lambda(ab)$ and involution $\Lambda(a)^\#:=\Lambda(a^*)$. We claim that this is a left Hilbert algebra, see, e.g.,~\cite{MR0526399}*{Chapter~10}. In order to see this we need to check two extra properties: $\UU^2$ is dense in~$\UU$ and the antilinear map $S_0\colon\UU\to H$, $\xi\mapsto\xi^\#$, is closable. The first property is immediate by condition~(c). In order to check the second property note that (iii) can be written~as
$$
(\Lambda(a^*),\Lambda(b))=(\Lambda(\sigma_{i\beta}(b^*)),\Lambda(a)).
$$
This shows that $\UU\subset D(S^*_0)$ and $S^*_0\Lambda(b)=\Lambda(\sigma_{i\beta}(b^*))$. Therefore $S^*_0$ is densely defined, so $S_0$ is closable. Denote by $S$ its closure. This is the modular involution associated to $\UU$.

It follows now that the von Neumann algebra $M:=\pi(A)''=\pi(\A)''\subset B(H)$ is in a standard form and $\UU$ defines a normal semifinite faithful weight $\Phi$ on $M$ such that $\Phi(\pi(a^*a))=\|\Lambda(a)\|^2=\varphi(a^*a)$ for all $a\in\A$. Therefore $\psi:=\Phi\circ\pi|_{A_+}$ is a lower semicontinuos densely defined weight satisfying $\psi(a^*a)=\varphi(a^*a)$ for all $a\in\A$. In particular, $\A\subset\NN_\psi$. By the polarization identity we then have $\psi(b^*a)=\varphi(b^*a)$ for all $a,b\in\A$, and using again property (c) we see that $\A\subset\M_\psi$ and $\psi(a)=\varphi(a)$ for all $a\in\A$.

In order to show that $\psi$ is a $\sigma$-KMS$_\beta$-weight we need to compute the modular group of $\Phi$. By property (i) we can define unitaries $u_t$ on $H$ by $u_t\Lambda(a):=\Lambda(\sigma_{-\beta t}(a))$ for $a\in\A$ and $t\in\R$. By property (ii) the one-parameter group $(u_t)_{t\in\R}$ is weakly operator continuous, hence it is also strongly operator continuous. Let $\Delta$ be its generator, so that $u_t=\Delta^{it}$.

Take $a\in\A$. We claim that the map $f_a\colon\C\to H$, $f_a(z):=\Lambda(\sigma_z(a))$, is analytic. From properties (i) and (ii) we see that it is at least continuous. Since the functions $z\mapsto (f_a(z),\Lambda(b))$ are analytic for all $b\in\A$ by (ii), analyticity of $f_a$ follows now by a standard argument. Namely, fix an open disk in $\C$ and let $\Gamma$ be its boundary. Then for all $z$ in the disk and all $b\in\A$ we have
$$
(f_a(z),\Lambda(b))=\frac{1}{2\pi i}\oint_\Gamma\frac{(f_a(w),\Lambda(b))}{w-z}dw.
$$
By density of $\Lambda(\A)$ in $H$ it follows that $f_a(z)=\frac{1}{2\pi i}\oint_\Gamma(w-z)^{-1}f_a(w)dw$, which implies analyticity.

We thus see that for every $a\in\A$ the vector $\Lambda(a)$ is analytic for $(u_t)_t$ and $\Delta^z\Lambda(a)=\Lambda(\sigma_{i\beta z}(a))$. In particular, $\Lambda(\A)$ is a core for any positive power of $\Delta$. Consider the antiunitary $J$ defined by $J\Lambda(a):=\Lambda(\sigma_{\frac{i\beta}{2}}(a^*))=\Lambda(\sigma_{-\frac{i\beta}{2}}(a)^*)$ for $a\in\A$. We then have $S\Lambda(a)=\Lambda(a^*)=J\Delta^{1/2}\Lambda(a)$. Hence $S=J\Delta^{1/2}$, so that $\Delta$ is the modular operator and $J$ is the modular conjugation associated to the Hilbert algebra $\UU$. Then for the modular group of $\Phi$ we have $\sigma_t^\Phi(\pi(a))=\Delta^{it}\pi(a)\Delta^{-it}=\pi(\sigma_{-\beta t}(a))$ on $\Lambda(\A)\subset H$ for all $a\in\A$, hence
$\sigma^\Phi_t\circ\pi=\pi\circ\sigma_{-\beta t}$ on $A$. Invariance of $\psi$ under~$\sigma$ for $\beta\ne0$ and the KMS-condition~\eqref{eq:KMS} follow now from the corresponding properties of~$\Phi$. For $\beta=0$ the invariance holds for similar reasons: the unitaries $v_t\Lambda(a):=\Lambda(\sigma_t(a))$ define automorphisms of $\UU$, hence the weight $\Phi$ is $\Ad v_t$-invariant, and since $(\Ad v_t)\circ\pi=\pi\circ\sigma_t$, we conclude that $\psi$ is $\sigma_t$-invariant.
\ep

Returning to groupoid C$^*$-algebras, it has been already shown by Christensen~\cite{MR4592883}*{Proposition~6.1} that $C_c(\G)\subset\M_\psi$ for any $\sigma^c$-KMS$_\beta$-weight $\psi$ on $C^*(\G)$. The following lemma and Proposition~\ref{prop:Ped} refine this result.

\begin{lemma}\label{lem:core}
Assume $\G$ is a locally compact \'etale groupoid such that the unit space~$\Gu$  is Hausdorff, $p\colon\LL\to\G$ is a Fell line bundle and $c\colon\G\to\R$ is a continuous $1$-cocycle. Then the subalgebra $\A:=\Gamma_c(\G,\LL)\subset C^*_r(\G,\LL)$ satisfies the conditions of Proposition~\ref{prop:Ped} with respect to the dynamics $\sigma^c$ on $C^*_r(\G,\LL)$ defined by~$c$.
\end{lemma}

\bp
If $W\subset\Gu$ is an open bisection and $\eta\in\Gamma_c(W,\LL)$, then $\eta$ is $\sigma^c$-analytic by~\eqref{eq:norm-equality2}, with $\sigma^c_z(\eta)=e^{izc}\eta\in\Gamma_c(W,\LL)$. This implies (a). Next, we have $C_c(\Gu)=\Ped(C_0(\Gu))\subset\Ped(C^*_r(\G,\LL)^{\sigma^c})$. For every $\eta\in\Gamma_c(\G,\LL)$, we can find $h\in C_c(\Gu)$ such that $\eta=h*\eta=\eta*h$. These two properties prove~(b) and~(c). Finally, (d) is satisfied by Lemma~\ref{lem:unibound2}.
\ep

In order to proceed we need to recall some terminology from measure theory on arbitrary Hausdorff locally compact spaces. Given such a space $X$, a Borel measure $\mu$ on $X$ is called a Radon measure if it is finite on every compact set, inner regular on every open set and outer regular on every Borel set. Then $\mu$ is inner regular on every Borel set of finite measure and, more generally, on the countable unions of such sets.

A subset $A\subset X$ is called \emph{locally Borel} if $A\cap K$ is Borel for every compact set $K\subset X$. Correspondingly, a function $f\colon X\to\C$ is called locally Borel if it is measurable with respect to the $\sigma$-algebra of locally Borel sets, equivalently, if $f|_K$ is Borel for every compact set $K\subset X$. A locally Borel set $A$ is said to be \emph{locally $\mu$-null}, or $\mu$-negligible, if $\mu(A\cap K)=0$ for every compact set $K\subset X$; clearly, then $\mu(A\cap Y)=0$ for every $\sigma$-compact subset $Y$ of $X$. One says that a property holds \emph{locally $\mu$-everywhere}, which we abbreviate as $\mu$-l.a.e., if the set of points $x\in X$ where it does not hold is locally $\mu$-null.

The key tool for deducing results for Radon measures from the case of finite measures is a \emph{concassage} of a Radon measure $\mu$. By definition, this is a collection $(K_\alpha)_\alpha$ of compact subsets of $X$ satisfying the following properties: (i) the sets~$K_\alpha$ are mutually disjoint; (ii) if an open subset $U\subset X$ intersects nontrivially some $K_\alpha$, then $\mu(U\cap K_\alpha)>0$; (iii) the set $X\setminus\bigcup_\alpha K_\alpha$ is locally $\mu$-null. Existence of a concassage is  proved by an application of Zorn's lemma, see~\cite{MR0426084}*{Theorem~I.13}; this argument will essentially be reproduced below in the setting of groupoids.

\smallskip

Assume now that $\G$ is a locally compact \'etale groupoid with Hausdorff unit space. We say that a Radon measure $\mu$ on $\Gu$ is quasi-invariant if $(T_W)_*(\mu|_{s(W)})\sim \mu|_{r(W)}$ for every open bisection $W\subset\G$ with relatively compact $r(W)$ and $s(W)$, where $T_W\colon s(W)\to r(W)$ is the homeomorphism defined by $W$. For $x\in r(W)$, we write $g^x_W$ for the unique element of $W\cap\G^x$. The following result shows that the Radon--Nikodym derivatives $\frac{d((T_W)_*\mu)}{d\mu}$ can be encoded by one function on $\G$.


\begin{prop}
Assume $\G$ is a locally compact \'etale groupoid with Hausdorff unit space~$\Gu$ and $\mu$ is a quasi-invariant Radon measure on $\Gu$. Then there exists a function $D\colon\G\to(0,+\infty)$ such that, for every open bisection $W\subset\G$ with relatively compact $r(W)$ and $s(W)$, the function~$D|_W$ is Borel and
$$
\frac{d((T_W)_*\mu)}{d\mu}(x)=D(g^x_W)^{-1}\ \ \text{for}\ \ \mu\text{-a.e.}\ x\in r(W).
$$
If $D'$ is another function with the same properties, then the Borel set $r(\{g\in W: D(g)\ne D'(g)\})$ is $\mu$-null for every $W$ as above.
\end{prop}

We will call any such function $D$ the Radon--Nikodym cocycle of~$\mu$.

\bp
The last statement of the proposition is clear, we only need to prove existence of~$D$. We claim that, similarly to existence of a concassage of a Radon measure, there exists a collection~$(K_\alpha)_\alpha$ of compact subsets of $\G$ satisfying the following properties:
\begin{enumerate}
  \item[(i)] each set $K_\alpha$ is contained in an open bisection;
  \item[(ii)] the sets $K_\alpha$ are mutually disjoint;
  \item[(iii)] if an open subset $U\subset\G$ intersects nontrivially some $K_\alpha$, then $\mu(r(U\cap K_\alpha))>0$;
  \item[(iv)] for every open bisection $W\subset\G$ with relatively compact $r(W)$, the set $r(W\setminus\bigcup_\alpha K_\alpha)\subset\Gu$ is $\mu$-null.
\end{enumerate}

In order to prove the claim, consider the set of all collections of compact sets satisfying conditions~(i)--(iii). The empty collection is contained in this set. By Zorn's lemma there is a maximal such collection $(K_\alpha)_\alpha$. We need to check that it satisfies (iv). Take an open bisection~$W$ such that $r(W)$ is relatively compact in~$\Gu$, so that $\mu(r(W))<\infty$. By condition (iii), there exist at most countably many indices~$\alpha_n$ such that $W\cap K_{\alpha_n}\ne\emptyset$. Hence $r(W\setminus\bigcup_\alpha K_\alpha)=r(W\setminus\bigcup_n K_{\alpha_n})$ is a Borel set of finite measure. If it is not of measure zero, then by regularity of~$\mu$ on Borel sets of finite measure there exists a compact subset $C\subset r(W\setminus\bigcup_\alpha K_\alpha)$ such that $\mu(C)>0$. Let $C_0\subset C$ be the support of $\mu|_C$ and put $K:=r^{-1}(C_0)\cap W$. Then the collection $(K_\alpha)_\alpha\cup\{K\}$ satisfies properties (i)--(iii), contradicting maximality of~$(K_\alpha)_\alpha$.

Denoting by $T_\alpha\colon s(K_\alpha)\to r(K_\alpha)$ the homeomorphisms defined by the bisections $K_\alpha$ and by~$g^x_\alpha$ the unique element of $K_\alpha\cap\G^x$ for $x\in r(K_\alpha)$, we put $D(g^x_\alpha):=\frac{d((T_\alpha)_*\mu)}{d\mu}(x)^{-1}$ for $x\in r(K_\alpha)$. For the remaining points $g\in\G\setminus\bigcup_\alpha K_\alpha$ we let $D(g):=1$. Since for every open bisection~$W$ with relatively compact $r(W)$ and $s(W)$ there are at most countably many indices $\alpha$ such that $W\cap K_\alpha\ne\emptyset$ and for every such $\alpha$ we have $\frac{d((T_W)_*\mu)}{d\mu}(x)=\frac{d((T_\alpha)_*\mu)}{d\mu}(x)$ for $\mu$-a.e.~$x\in r(W\cap K_\alpha)$, we see that $D$ has the required properties.
\ep

\begin{remark}
The chain rule for Radon--Nikodym derivatives implies that if $W$ is an open bisection with relatively compact $r(W)$ and $s(W)$, then the set
$r(\{g\in W: D(g)^{-1}\ne D(g^{-1})\})$ is $\mu$-null, and if $V$ is another such bisection, then the set
$$
r(\{gh\in VW: g\in V,\ h\in W,\ s(g)=r(h),\ D(gh)\ne D(g)D(h)\})
$$
is $\mu$-null. Replacing $D(g)$ by $D(g)^{1/2}D(g^{-1})^{-1/2}$, we can in fact get the identity $D(g)^{-1}=D(g^{-1})$ to hold for all~$g\in\G$. \ee
\end{remark}

Assume we are given a Fell line bundle $p\colon\LL\to \G$ and, for every $x\in\Gu$, a state $\varphi_x$ on the twisted group C$^*$-algebra $C^*(\Gxx,\LL)$.

\begin{defn}
We say that the family $(\varphi_x)_x$ is \emph{locally Borel} if the function $\Gu\ni x\mapsto\sum_{g\in\Gxx}\varphi_x(\eta(g))$ is Borel for all $\eta\in\Gamma_c(\G,\LL)$.
\end{defn}

We are now ready to prove that following generalization of  \cite{MR3138368}*{Theorem~1.3}.

\begin{thm}\label{thm:KMS}
Assume $\G$ is a locally compact \'etale groupoid with Hausdorff unit space~$\Gu$ and~$p\colon\LL\to\G$ is a Fell line bundle such that $\G$ can be covered by countably many open bisections over which $\LL$ admits unitary sections. Assume $c\colon\G\to\R$ is a continuous $1$-cocycle and consider the corresponding dynamics $\sigma^c$ defined by~\eqref{eq:dynamics}. Then, for every $\beta\in \R$, there exists a one-to-one correspondence between the $\sigma^c$-KMS$_\beta$-weights on~$C^*(\G,\LL)$ and the equivalence classes of pairs $(\mu,(\tau_x)_{x\in \Gu})$ consisting of a Radon measure $\mu$ on $\Gu$ and a locally Borel family of tracial states~$\tau_x$ on $C^*(\G^x_x,\LL)$ such that
\begin{enumerate}
\item[(i)] $\mu$ is quasi-invariant, with Radon--Nikodym cocycle $e^{-\beta c}$;
\item[(ii)] $\tau_x=\tau_{r(g)}\circ\Ad u$ for $\mu$-l.a.e.~$x$, all $g\in \G_x$ and  some (equivalently, every) unitary element $u\in \LL_g$;
\item[(iii)] $\tau_x|_{\LL_g}=0$ for $\mu$-l.a.e.~$x$ and all $g\in \G^x_x\setminus c^{-1}(0)$.
\end{enumerate}
Specifically, the weight corresponding to $(\mu,(\tau_x)_x)$ is given~by
\begin{equation}\label{eq:KMS-weight}
\psi(\eta)=\int_{\Gu}\sum_{g\in \Gxx}\tau_x(\eta(g))d\mu(x)\ \ \hbox{for}\ \ \eta\in \Gamma_c(\G,\LL).
\end{equation}
Here we say that two pairs $(\mu,(\tau_x)_{x})$ and $(\mu',(\tau'_x)_{x})$ are equivalent if $\mu=\mu'$ and $\tau_x=\tau'_x$ $\mu$-l.a.e.
\end{thm}

We remark that condition (iii) is redundant when $\beta\ne0$. In order to see this, take an open bisection $W$ with relatively compact $r(W)$ and $s(W)$. As before, consider the corresponding homeomorphism $T_W\colon s(W)\to r(W)$ and let $g^x_W$ be the unique element of $W\cap\G^x$ for $x\in r(W)$. Consider the locally closed set $X_W$ of points $x\in s(W)$ fixed by $T_W$. Then, on the one hand,  $\frac{d\left((T_{W})_*\mu\right)}{d\mu}(x)=1$ for $\mu$-a.e.~$x\in X_W$. On the other hand, by condition (i), we have  $\frac{d\left((T_{W})_*\mu\right)}{d\mu}(x)=e^{\beta c(g^x_W)}$ for $\mu$-a.e.~$x\in X_W$. Therefore $c(g^x_W)=0$ for $\mu$-a.e.~$x\in X_W$, or in other words, the set $\{x\in\Gu: W\cap(\Gxx\setminus c^{-1}(0))\ne\emptyset\}$ is $\mu$-null. By regularity of $\mu$ on Borel sets of finite measure this implies that for any open bisection $W$ the Borel set $\{x\in\Gu: W\cap(\Gxx\setminus c^{-1}(0))\ne\emptyset\}$ is locally $\mu$-null. Since $\G$ can be covered by countably many open bisections, it follows that the set $\{x\in\Gu: \Gxx\not\subset c^{-1}(0)\}$ is locally $\mu$-null, so condition (iii) is a trivial consequence of (i).

\bp[Proof of Theorem~\ref{thm:KMS}]
The proof is similar to that of \cite{MR3138368}*{Theorem~1.3}, see also \cite{MR4222432}*{Theorem~3.4} and~\cite{MR4592883}*{Theorem~6.3}, so instead of giving full details we will concentrate on the parts that require some changes.

As a first step let us show that there is a one-to-one correspondence, given by~\eqref{eq:KMS-weight}, between the positive linear functionals on $\Gamma_c(\G,\LL)$ such that $C_c(\Gu)$ is contained in the centralizer of~$\psi$, meaning that $\psi(f*\eta)=\psi(\eta*f)$ for all $f\in C_c(\Gu)$ and $\eta\in\Gamma_c(\G,\LL)$, and the equivalence classes of pairs $(\mu,(\tau_x)_x)$ consisting of a Radon measure $\mu$ and a locally Borel family of states~$\tau_x$ on $C^*(\Gxx,\LL)$. This is a generalization of \cite{MR3138368}*{Theorem~1.1}.

Assume first that $\psi$ is a positive linear functional on $\Gamma_c(\G,\LL)$ such that $C_c(\Gu)$ is contained in the centralizer of $\psi$. The restriction of $\psi$ to $C_c(\Gu)$ defines a Radon measure $\mu$ on $\Gu$. Consider the GNS-map $\Lambda\colon \Gamma_c(\G,\LL)\to H$ associated to $\psi$ and the corresponding representation $\pi\colon C^*(\G,\LL)\to B(H)$. Fix a relatively compact open subset $U\subset\Gu$ and a function $h\in C_c(\Gu)$ such that $\un_U\le h\le1$. Then the functional $\psi_h(\eta):=\psi(h*\eta *h)$ is still positive and contains $C_c(\Gu)$ in its centralizer. We have $\psi_h=(\pi(\cdot)\Lambda(h),\Lambda(h))$, so $\psi_h$ extends to a positive linear functional on $C^*(\G,\LL)$. Applying Theorem~\ref{thm:main-twist} to the subrepresentation of $\pi$ on $\overline{\pi(C^*(\G,\LL))\Lambda(h)}$ we get a quasi-invariant regular Borel probability measure $\nu_h$ on $\Gu$ and a nondegenerate Borel action $\LL\curvearrowright\HH$ such that $\Lambda(h)$ corresponds to a square-integrable section $(\xi_x)_x$ of $\HH$ and
\begin{equation}\label{eq:cent1}
\psi_h(\eta)=\int_{\Gu}\sum_{g\in \G_x}D_h(g)^{-1/2}(L(\eta(g))\xi_{s(g)},\xi_x)d\nu_h(x)\ \ \hbox{for}\ \ \eta\in \Gamma_c(\G,\LL),
\end{equation}
where $D_h$ is the Radon--Nikodym cocycle of $\nu_h$.

Take an open bisection $W$ over which there is a unitary section $\eta'\colon W\to\LL$. For every point $x_0\in r(W)$ not fixed by $T_W^{-1}$, choose a relatively compact open neighborhood $V_{x_0}\subset r(W)$ such that $V_{x_0}\cap T_W^{-1}(V_{x_0})=\emptyset$. Then for every $f\in C_c(V_{x_0})$ we have $\eta'*f=0$. It follows that $\psi_h(f*\eta')=0$. This implies that $(L(\eta'(g^x_W))\xi_{T_W^{-1}(x)},\xi_x)=0$ for $\nu_h$-a.e.~$x\in V_{x_0}$. By regularity of the measure $\nu_h$ we conclude that the same is true for $\nu_h$-a.e.~$x\in r(W)$ not fixed by $T_W^{-1}$. Note also that for the same reason as that (iii) follows from (i) for $\beta\ne0$, we have $D_h(g^x_W)^{-1}=\frac{d\left((T_{W})_*\nu_h\right)}{d\nu_h}(x)=1$ for $\nu_h$-a.e.~$x$ fixed by $T_W$. It follows that~\eqref{eq:cent1} can be written as
\begin{equation*}\label{eq:cent2}
\psi_h(\eta)=\int_{\Gu}\sum_{g\in \Gxx}(L(\eta(g))\xi_x,\xi_x)d\nu_h(x)\ \ \hbox{for}\ \ \eta\in \Gamma_c(\G,\LL).
\end{equation*}

By restricting to $C_c(\Gu)$ we see that we must have the equality of measures $\|\xi_x\|^2d\nu_h(x)=h(x)^2d\mu(x)$. It follows that for $\mu$-a.e.~$x$ with $h(x)>0$ we get a well-defined state $\tau_x:=\|\xi_x\|^{-2}(L(\cdot)\xi_x,\xi_x)$ on $C^*(\Gxx,\LL)$, and then
\begin{equation*}
\psi_h(\eta)=\int_{\Gu}\sum_{g\in \Gxx}\tau_x(\eta(g))h(x)^2d\mu(x)\ \ \hbox{for}\ \ \eta\in \Gamma_c(\G,\LL).
\end{equation*}
This shows that with our choice of $\tau_x$ for $x\in U$ identity~\eqref{eq:KMS} holds for all sections $\eta\in \Gamma_c(\G,\LL)$ that are concentrated on $\G|_U:=r^{-1}(U)\cap s^{-1}(U)$. It is easy to show that the collection $(\tau_x)_{x\in U}$ is essentially uniquely determined by this.

We thus see that it is possible to choose the required states $\tau_x$ for $\mu$-a.e.~$x$ in any given relatively compact open set $U$, and the choice is essentially unique, so that if $(\tau_x)_{x\in U}$ are $(\tau_x')_{x\in V}$ are our choices for $U$ and $V$, then $\tau_x=\tau_x'$ for $\mu$-a.e.~$x\in U\cap V$. It remains to explain how to make a global choice that is locally Borel. For this take a concassage $(K_\alpha)_\alpha$ of the Radon measure~$\mu$. Every compact set $K_\alpha$ is contained in some relatively compact open set $U_\alpha$, hence we can choose Borel families of states $\tau_x$ for $x\in K_\alpha\setminus N_\alpha$, where $N_\alpha\subset K_\alpha$ are certain Borel sets of measure zero, by first choosing such families for $U_\alpha$. For the remaining points $x\in\Gu\setminus\bigcup_\alpha(K_\alpha\setminus N_\alpha)$ we take the canonical trace on $C^*(\Gxx,\LL)$, that is, $\tau_x(\eta)=\eta(x)$ for $\eta\in\Gamma_c(\Gxx,\LL)$, where we as usual identify $\LL_x$ with $\C$.

In the opposite direction we need to show that \eqref{eq:KMS-weight} defines a positive linear functional on $\Gamma_c(\G,\LL)$ with centralizer containing $C_c(\Gu)$ for any pair $(\mu,(\tau_x)_x)$ consisting of a Radon measure~$\mu$ and a locally Borel family of states~$\tau_x$. For this it suffices to consider delta-measures $\mu=\delta_x$, in which case we need to show that any state $\tau$ on $C^*(\Gxx,\LL)$ defines a positive linear functional $\psi$ on $\Gamma_c(\G,\LL)$ such that $\psi(\eta)=\sum_{g\in\Gxx}\tau(\eta(g))$. This is proved by inducing the GNS-representation of $C^*(\Gxx,\LL)$ associated to $\tau$ to a representation of~$C^*(\G,\LL)$, see \cites{MR3138368,MR4222432} for details. That $C_c(\Gu)$ is contained in the centralizer of $\psi$ is obvious. This finishes the first step of the proof.

\smallskip

Now, in order to prove the theorem, by Proposition~\ref{prop:Ped} and Lemma~\ref{lem:core} we need to show that a positive linear functional $\psi$ on $\A:=\Gamma_c(\G,\LL)$ of the form~\eqref{eq:KMS-weight} satisfies conditions (i)--(iii) of Proposition~\ref{prop:Ped} if and only if conditions (i)--(iii) of the theorem are satisfied. It is not difficult to see that condition (ii) of Proposition~\ref{prop:Ped} is always satisfied for such $\psi$, while condition (i) there is equivalent to condition (iii) of the theorem. It is checked then that conditions~(i) and~(ii) of the theorem are equivalent to condition (iii) of Proposition~\ref{prop:Ped}, see again \cites{MR3138368,MR4222432} for details.
\ep

It is clear that this theorem gives  also a description of $\sigma^c$-KMS$_\beta$-states on $C^*(\G,\LL)$ - they correspond to the pairs $(\mu,(\tau_x)_x)$ such that $\mu$ is a probability measure.

\smallskip

The theorem can in principle be used to study KMS-weights and states for arbitrary $\G$. As a sample application let us prove the following generalization of \cite{MR4411323}*{Corollary~2.4}. For tracial weights the same result has been recently obtained in~\cite{MS}*{Section~4}.

\begin{cor}
Assume $\G$ is a locally compact \'etale groupoid with Hausdorff unit space~$\Gu$, $p\colon\LL\to\G$ is a Fell line bundle, $c\colon\G\to\R$ is a continuous $1$-cocycle, $\beta\in\R$ and $\mu$ is a quasi-invariant Radon measure on $\Gu$ with Radon--Nikodym cocycle $e^{-\beta c}$. Assume also that either $\beta\ne0$ or $c=0$. Consider the following conditions:
\begin{itemize}
  \item[(1)] there is a unique $\sigma^c$-KMS$_\beta$-weight on $C^*(\G,\LL)$ such that its restriction to $C_c(\Gu)$ is defined by $\mu$;
  \item[(2)] for every open bisection $W\subset\G$, the Borel set $\{x\in \Gu: W\cap(\Gxx\setminus\{x\})\ne\emptyset\}$ is locally $\mu$-null.
\end{itemize}
Then $(2)\Rightarrow(1)$. If the Fell bundle $\LL$ is trivial, then the two conditions are equivalent.
\end{cor}

\bp
Choose a directed family of open subgroupoids $\G(i)\subset\G$ ($i\in I$) such that $\bigcup_i\G(i)=\G$ and each~$\G(i)$  can be covered by countably many open bisections over which $\LL$ admits unitary sections. Replacing $\G(i)$ by $\G(i)\cup\Gu$ we may assume that $\G(i)^{(0)}=\Gu$ for all~$i\in I$. As the algebra $\Gamma_c(\G,\LL)$ is the increasing union of the subalgebras $\Gamma_c(\G(i),\LL)$, by Proposition~\ref{prop:Ped} and Lemma~\ref{lem:core} applied to the groupoids $\G$ and $\G(i)$ we get a one-to-one correspondence between the $\sigma^c$-KMS$_\beta$-weights $\psi$ on $C^*(\G,\LL)$ and the $I$-indexed collections of $\sigma^c$-KMS$_\beta$-weights $\psi_i$ on $C^*(\G(i),\LL)$ such that $\psi_i=\psi_j$ on $\Gamma_c(\G(i),\LL)$ whenever $\G(i)\subset\G(j)$. By Theorem~\ref{thm:KMS} applied to the groupoids~$\G(i)$ we conclude then that there is a one-to-one correspondence between the $\sigma^c$-KMS$_\beta$-weights $\psi$ on $C^*(\G,\LL)$ such that $\psi|_{C_c(\Gu)}$ is defined by $\mu$ and the equivalence classes of $I$-indexed collections of locally Borel families $(\tau^i_x)_{x\in\Gu}$ of tracial states $\tau^i_x$ on $C^*(\G(i)^x_x,\LL)$ satisfying conditions (ii) and (iii) of the theorem (for $\G(i)$ in place of $\G$) such that $\tau^i_x=\tau^j_x|_{C^*(\G(i)^x_x,\LL)}$ for $\mu$-l.a.e.~$x$ whenever $\G(i)\subset\G(j)$. Here we say that two $I$-indexed collections of families $(\tau^i_x)_{x\in\Gu}$ and $(\tilde\tau^i_x)_{x\in\Gu}$ are equivalent if for every $i\in I$ we have $\tau^i_x=\tilde\tau^i_x$ for $\mu$-l.a.e.~$x$.

The canonical traces $\tau^i_x$ on $C^*(\G(i)^x_x,\LL)$ obviously form such a collection of families. If condition (2) is satisfied, then for every $i\in I$ the set $\{x\in\Gu:\G(i)^x_x\ne\{x\}\}$ is locally $\mu$-null, so up to equivalence there are no other such collections. Thus, (2)$\Rightarrow$(1).

Assume now that $\LL$ is trivial. Then the trivial characters $\eps^i_x$ on $\G(i)^x_x$ also form such a collection of families. Note that here we use that condition (iii) in Theorem~\ref{thm:KMS} is trivially satisfied  for $\beta\ne0$, while if $\beta=0$, then $c=0$ by assumption. Therefore if (1) holds, then for every $i$ we must have that  $\tau^i_x=\eps^i_x$ for $\mu$-l.a.e.~$x$, that is, $\G(i)^x_x=\{x\}$ for $\mu$-l.a.e.~$x$. Equivalently, for every $i$ and every open bisection $W\subset\G(i)$, the set $\{x\in \Gu: W\cap(\G(i)^x_x\setminus\{x\})\ne\emptyset\}$ is locally $\mu$-null. Since every compact subset of $\G$ is contained in some $\G(i)$, the last condition is equivalent to (2).
\ep

\bigskip

\end{document}